\documentclass[a4paper,12pt]{article}
\usepackage{amsmath,amsthm,amssymb,amsfonts}
\usepackage{bbm,bm}
\usepackage{latexsym}
\usepackage{mathrsfs}
\usepackage{threeparttable}
\usepackage{tabularx}
\usepackage{booktabs}
\usepackage{graphicx,subfigure}
\usepackage{color}
\usepackage{indentfirst}
\usepackage{geometry}
\usepackage{caption}
\usepackage{lineno}
\usepackage{abstract}
\usepackage{enumerate,mdwlist}
\usepackage[numbers,sort&compress]{natbib}
\usepackage[colorlinks,
            linkcolor=blue, 
            citecolor=red  
            ]{hyperref}
\usepackage{epstopdf}

\def \[{\begin{equation}}
\def \]{\end{equation}}

\newtheorem{thm}{Theorem}[section]

\newtheorem{lem}[thm]{Lemma}

\newtheorem{conj}[thm]{Conjecture}

\begin{document}

\setlength{\baselineskip}{20pt}
\begin{center}{\Large \bf Minimum degree of minimal (\emph{n}-10)-factor-critical graphs}\footnote{This work
is supported by NSFC\,(Grant No.  12271229).}

\vspace{4mm}

{Jing Guo, Heping Zhang \footnote{The corresponding author.}
\renewcommand\thefootnote{}\footnote{E-mail addresses: guoj20@lzu.edu.cn (J. Guo), zhanghp@lzu.edu.cn (H. Zhang).}}

\vspace{2mm}

\footnotesize{ School of Mathematics and Statistics, Lanzhou University, Lanzhou, Gansu 730000, P. R. China}

\end{center}
\noindent {\bf Abstract}: A graph $G$ of order $n$ is said to be $k$-factor-critical
 for integers $1\leq k < n$, if the removal of any $k$
 vertices results in a graph with a perfect matching.
 A $k$-factor-critical graph $G$
 is called minimal if for any edge $e\in E(G)$, $G-e$ is not $k$-factor-critical.
 In 1998, O. Favaron and M. Shi conjectured that every minimal $k$-factor-critical graph
 of order $n$ has the minimum degree $k+1$ and confirmed it for $k=1, n-2, n-4$ and $n-6$.
 By using a novel approach, we have confirmed it for $k = n - 8$ in a previous paper.
 Continuing this method, we prove  the conjecture to be true for $k=n-10$ in this paper.

\vspace{2mm} \noindent{\bf Keywords}: Perfect matching;
Minimal $k$-factor-critical graph; Minimum degree.
\vspace{2mm}

\noindent{AMS subject classification:} 05C70,\ 05C07

 {\setcounter{section}{0}
\section{Introduction}\setcounter{equation}{0}

Only finite and simple graphs are considered in this article. Let $G$ be a
graph with vertex set $V(G)$ and edge set $E(G)$.
The {\em order} of $G$ is the cardinality of $V(G)$.
A set of edges $M\subseteq E(G)$ is called a {\em matching} of $G$ if no two of them share an end-vertex.
A matching of $G$ is said to be a {\em perfect matching} or a 1{\em -factor} if it covers all vertices of $G$.
The concepts of factor-critical and bicritical graphs were introduced by T. Gallai \cite{GT} and
L. Lov\'{a}sz \cite{LL}, respectively.
A graph $G$ is called {\em factor-critical} if the removal of any vertex
of $G$ results in a graph with a perfect matching.
A graph $G$ with at least one edge is called {\em bicritical} if the removal
of any pair of distinct vertices of $G$ results in a graph with a perfect matching.

A 3-connected bicritical graph is the so-called {\em brick},
which plays a key role in matching theory of graphs.
J. Edmonds et al. \cite{ELW} and L. Lov\'{a}sz \cite{LO} proposed and developed the
``tight set decomposition" of matching-covered graphs into list of bricks in an essentially unique manner.
The decomposition can reduce some matching problems of graphs to bricks,
such as,  the dimension of matching lattices \cite{LO} and  perfect matching polytopes \cite{ELW},
 Pfaffian orientation \cite{LR, VY}, etc.

Generally, O. Favaron \cite{F} and Q. Yu \cite{Y} introduced, independently,
$k$-factor-critical graphs as a generalization of factor-critical and bicritical graphs.
A graph $G$ of order $n$ is said to be {\em $k$-factor-critical} for integers $1\leq k < n$,
if the removal of any $k$ vertices results in a graph with a perfect matching.
They gave characterizations of $k$-factor-critical graphs and
the following important property on connectivity.

\begin{thm}[\cite{F, Y}]\label{conn}
If $G$ is $k$-factor-critical for some $1\leq k< n$ with $n+k$ even,
then $G$ is $k$-connected, $(k+1)$-edge-connected and $(k-2)$-factor-critical if $k\geq 2$.
\end{thm}

For more about the $k$-factor-critical graphs, the reader is referred to
articles \cite{PS, PMD, N, LDYQ, ZWZ} and a monograph \cite{YL}.


A graph $G$ is called {\em minimal $k$-factor-critical} if $G$ is $k$-factor-critical
but $G-e$ is not  $k$-factor-critical for any $e\in E(G)$.
L. Lov\'{a}sz and M. D. Plummer \cite{LP, LP75} considered minimal bicritical graphs and revealed some excluded subgraphs (wheel and $K_{3,3}$). For a graph $G$ with a vertex $v$, let  $d_G(v)$ denote the degree of $v$ in $G$, the number of edges incident with vertex $v$, and  $\delta(G)$ the {\em minimum degree} of $G$.
O. Favaron and M. Shi \cite{FS} studied the minimum degree of
minimal $k$-factor-critical graphs and obtained the following result.

\begin{thm}[\cite{FS}]\label{Special}
Let $G$ be a minimal $k$-factor-critical graph of order $n$.
If $k=1, n-2, n-4$ or $n-6$, then $\delta(G)=k+1$.
\end{thm}






Since every $k$-factor-critical graph is $(k+1)$-edge-connected,
it has minimum degree at least $k+1$. So in 1998 they posed a problem:
does Theorem \ref{Special} hold for general $k$?
Afterward, Z. Zhang et al. \cite{ZWL} formally proposed the following conjecture.


\begin{conj}[\cite{FS,ZWL}]\label{Conj}
Let $G$ be a minimal $k$-factor-critical graph of order $n$  with $1\leq k<n$.
Then $\delta(G)=k+1$.
\end{conj}

A closely related concept to $k$-factor-critical is that of $q$-extendable.
D. Lou and Q. Yu \cite{LY} conjectured that
any minimal $q$-extendable graph $G$ on $n$ vertices with $n\leq 4q$
has minimum degree $q + 1, 2q$ or $2q + 1$.
Z. Zhang et al. \cite{ZWL} pointed out that the conjecture  is actually a part of Conjecture \ref{Conj}
except the case $n=4q$.

A brick $G$ is {\em minimal} if $G-e$ is not a brick for every edge $e$ of $G$.
In 1973, L. Lov\'{a}sz early conjectured that every minimal brick has two adjacent
vertices of degree three.
S. Norine and R. Thomas \cite{NT} presented a recursive procedure for
generating minimal bricks and  obtained that
every minimal brick has at least three vertices of degree three.  Further, F. Lin et al. [14] showed that
every minimal brick has at least four vertices of degree three.
For other results on minimal bricks, we refer to \cite{CLM, LZL, BS}. From such results we have that a 3-connected minimal bicritical graph has the minimum degree three since it is also a minimal brick.  However Conjecture \ref{Conj} remains open even for $k=2$.

In a previous paper \cite{GZ},  we   considered Conjecture \ref{Conj} for  large $k$.
By using a novel  method we confirmed  Conjecture \ref{Conj} to be true  not only for $k = n - 4$ and $n-6$ but also for $k=n-8$.
Continuing this method, in this article we confirm that Conjecture \ref{Conj} holds for $k=n-10$ and obtain our main theorem as follows.

\begin{thm}[Main Theorem]\label{Main}
If $G$ is minimal $(n-10)$-factor-critical graph of order $n\geq 12$,
then $\delta(G)=n-9$.
\end{thm}

In next section  some preliminaries are given. Section 3 is devoted to
a detailed proof of Theorem \ref{Main}.

\section{Some preliminaries}

For any set $X\subseteq V(G)$, $G[X]$ denotes the subgraph of $G$ induced by $X$,
and $G-X=G[V(G)-X]$.
For an edge $e=uv\in E(G)$, $G-e$ or $G-uv$ stands for the graph $(V(G), E(G)-\{e\})$.
Similarly, if $u, v\in V(G)$ are nonadjacent vertices of $G$,
$G+uv$ stands for the graph $(V(G), E(G)\cup \{e\})$.
A vertex of $G$ with degree one is called a {\em pendent vertex}.
An {\em independent set} of a graph is a set of pairwise nonadjacent vertices.
The {\em complete graph} $K_{n}$ is the graph of order $n$ in which any two vertices are adjacent. A graph is {\em nontrivial} if it has order at least two.

The following is Tutte's $1$-factor Theorem. As usual we let $C_{o}(G)$ denote the number of
odd components of a graph $G$.

\begin{thm}[\cite{TTW}]\label{Tutte}
A graph $G$ has a $1$-factor
if and only if $C_{o}(G-X)\leq |X|$ for any $X \subseteq V(G)$.
\end{thm}

A  stronger result was presented in \cite{LP}
which we make use of in our proof.

\begin{thm}[\cite{LP, TTW}]\label{tutte} A graph $G$ has no $1$-factor
if and only if there exists $X \subseteq V(G)$ such that all components of $G-X$ are
factor-critical and $C_{o}(G-X)\geq |X|+2$.
\end{thm}

The property of $k$-factor-critical graphs is presented as follows,
which were obtained by O. Favaron \cite{F} and Q. Yu \cite{Y}, independently.

\begin{lem}[\cite{F, Y}]\label{FY}
A graph $G$ is $k$-factor-critical if and only if $C_{o}(G-B)\leq |B|-k$ for
any $B \subseteq V(G)$ with $|B|\geq k$.
\end{lem}

O. Favaron and M. Shi \cite{FS} characterized minimal $k$-factor-critical graphs.

\begin{lem}[\cite{FS}]\label{minimal}
Let $G$ be a $k$-factor-critical graph. Then $G$ is minimal
if and only if for each $e=uv\in E(G)$, there exists $S_{e}\subseteq V(G)-\{u,v\}$ with $|S_{e}|=k$
such that every perfect matching of $G-S_{e}$ contains $e$.
\end{lem}

For a graph, the {\em neighborhood} of a vertex $x$ is $N(x):=\{y \mid y\in V(G),$ $xy\in E(G)\}$,
and the {\em closed neighborhood} is $N[x]:=N(x)\cup \{x\}$.
Then $\overline{N[x]}:=V(G)\setminus N[x]$ is called the {\em non-neighborhood} of $x$ in $G$,
which will play a critical role in subsequent discussions.

\section{Proof of Theorem \ref{Main}}

We first give a sketch for the lengthy proof of Theorem \ref{Main}. We proceed  by contradiction.
Since $G$ is minimal $(n-10)$-factor-critical graph, for every edge $e\in E(G)$,
$G-e$ is not $(n-10)$-factor-critical. By Lemma \ref{minimal} there exists a set $S_{e}\subseteq V(G)$
with $|S_{e}|=n-10$ such that $G_e=G-e-S_{e}$ has no perfect matchings.
By the stronger Tutte's $1$-factor Theorem, we have total fourteen configurations of $G-e-S_{e}$ which has order 10.
By analysing some properties of common non-neighborhood of the end-vertices of an edge, for each configuration  we always find a suitable (other) edge $e'$ so that
$G-e'-S_{e'}$ is not any one of the fourteen configurations, which yields a   contradiction.

We are now ready to prove our main theorem.

\smallskip
\noindent{\bf Proof of Theorem \ref{Main}.}
By Lemma \ref{conn}, $\delta(G)\geq n-9$.
Suppose to the contrary that $\delta(G)\geq n-8$.

\smallskip
{\textbf{Claim 1.}
For every $e=uv\in E(G)$, there exists $S_{e}\subseteq V(G)-\{u, v\}$
with $|S_{e}|=n-10$ such that $G_{e}=G-e-S_{e}$ has no perfect matchings. Further, $G_{e}$
is one of Configurations $C1$ to $C14$ (relative to edge $e$) as shown in Fig. 1.
(We bear in mind that notations $S_{e}$ and $G_{e}$ always are used in such meanings in next discussions.)}

\begin{figure}[h]
\centering
\includegraphics[height=9cm,width=15cm]{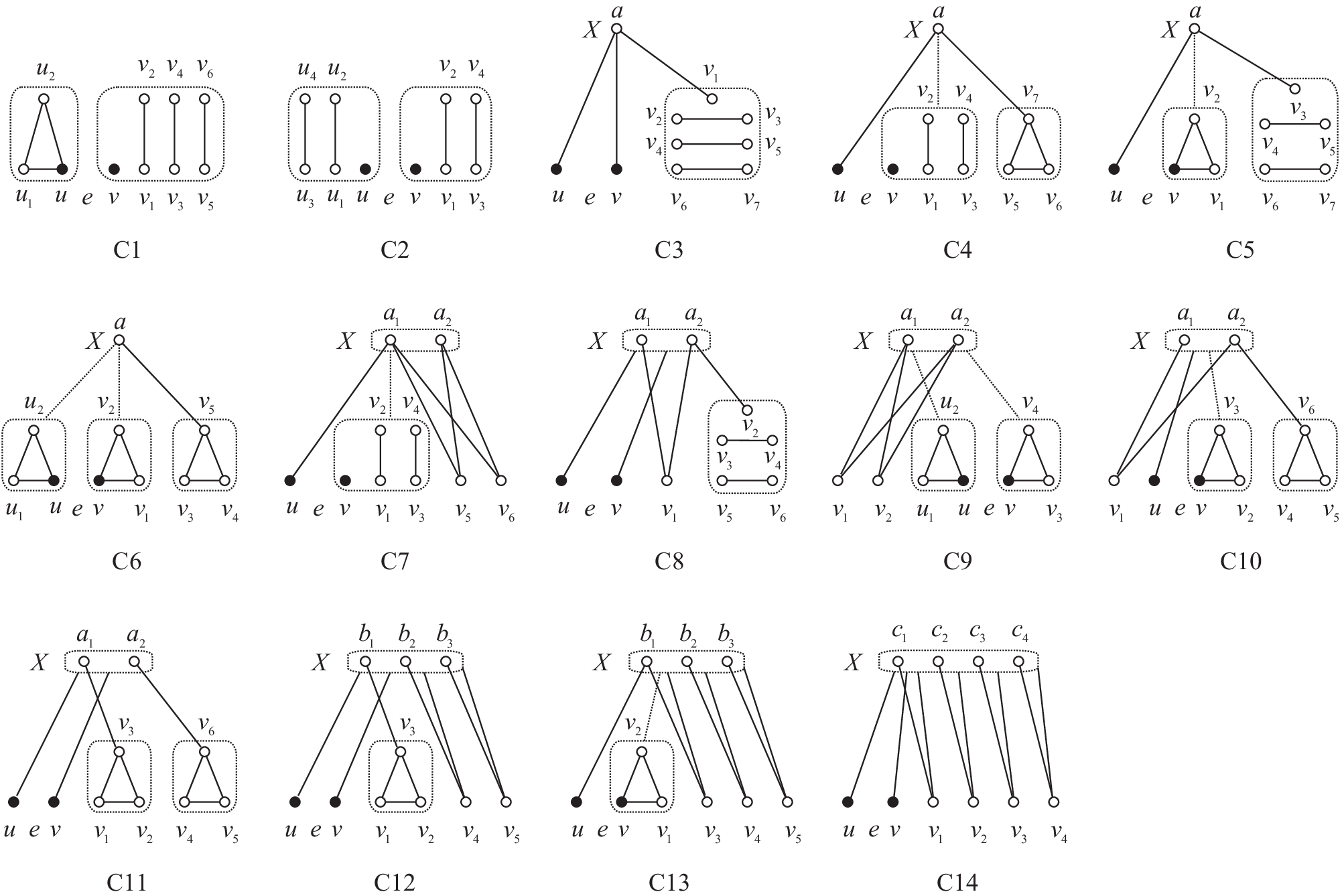}
\caption{\label{tu-01} The fourteen configurations of $G_{e}=G-e-S_{e}$. }
(The vertices within a dotted box induce a factor-critical subgraph
and dotted edge indicates an optional edge.)
\end{figure}

Since $G$ is minimal $(n-10)$-factor-critical graph, by Lemma \ref{minimal},
for any $e=uv\in E(G)$, there exists $S_{e}\subseteq V(G)-\{u, v\}$ with $|S_{e}|=n-10$
such that every perfect matching of $G-S_{e}$ contains $e$. Let $G_{e}=G-e-S_{e}$.
Then $G_{e}$ has order $10$ and no perfect matchings.
By Theorem \ref{tutte}, there exists $X\subseteq V(G_{e})$ such that
all components of $G_{e}-X$ are factor-critical and $C_{o}(G_{e}-X)\geq |X|+2$.
So $|X|+2\leq$ $C_{o}(G_{e}-X)\leq$ $|V(G_{e}-X)|$$=10-|X|$. Thus $|X|\leq 4$.
Since $G_{e}+e=G-S_{e}$ has a $1$-factor, $C_{o}(G_{e}-X)=|X|+2$ and $u$ and $v$ belong respectively to
two distinct odd components of $G_{e}-X$. Moreover, $\delta(G-S_{e})\geq 2$.
Then $G_{e}+e=G-S_{e}$ has no pendent vertex. So $G_{e}$ has no isolated vertex.

If $|X|=0$, then $G_{e}$ has exactly two odd components.
Since each component of $G_{e}$ is a factor-critical graph with at least three vertices,
$G_{e}$ has two possible cases as configurations $C1$ and $C2$.

If $|X|=1$, then $C_{o}(G_{e}-X)=3$ and $G_{e}-X$ has at most two trivial odd components.
If $G_{e}-X$ has two trivial odd components, then $e$ must join them.
Otherwise, $G_{e}+e$ has a pendent vertex, a contradiction.
Further, the third component is a factor-critical graph with seven vertices, so $G_{e}$ is $C3$.
If $G_{e}-X$ has only one trivial odd component, then the other two
odd components have three and five vertices, respectively. Thus
$e$ joins the trivial odd component and a nontrivial odd component. So $G_{e}$ is $C4$ or $C5$.
If $G_{e}-X$ has no trivial odd component, then each of the three odd components has three vertices.
So $G_{e}$ is $C6$.

If $|X|=2$, then $C_{o}(G_{e}-X)=4$ and $G_{e}-X$ has two or three trivial odd components.
If $G_{e}-X$ has three trivial odd components, then the other has
five vertices. So $G_{e}$ is $C7$ or $C8$ according to the possible position of edge $e$.
If $G_{e}-X$ has two trivial odd components, then each of the others
is a $K_{3}$, so $G_{e}$ is $C9$, $C10$ or $C11$.

If $|X|=3$, then $C_{o}(G_{e}-X)=5$. Thus $G_{e}-X$ consists of four trivial odd components and
a $K_{3}$. So $G_{e}$ is $C12$ or $C13$.

If $|X|=4$, then $C_{o}(G_{e}-X)=6$ and $G_{e}-X$ consists of six trivial odd components. So $G_{e}$ is $C14$.
Thus Claim 1 holds.

\smallskip

For every $x\in V(G)$, $\overline{N[x]}$ has at most seven vertices in $V(G)$ as $d_{G}(x)\geq n-8$.
For each configuration discussed below, let $C_{x}$ denote the odd component of $Ci-X$ containing vertex $x$
for $i=1, 2, \ldots, 14$. Note that since every $C_{x}$ is factor-critical,
any vertex of $C_{x}$ has at least two neighbors in $C_{x}$ unless $C_{x}$ is trivial.

\smallskip
{\bf Claim 2.} The non-neighborhoods of $u$ and $v$ have the following possible intersections:

\indent$(1)$ If $G_{e}$ is $C1$ or $C2$, then $|\overline{N[u]}\cap \overline{N[v]}|\leq 5$;

\indent$(2)$ If $G_{e}$ is $C3$, then $|\overline{N[u]}\cap \overline{N[v]}|=7$;

\indent$(3)$ If $G_{e}$ is $C4$, $C6$ or $C13$, then $3 \leq |\overline{N[u]}\cap \overline{N[v]}|\leq 5$;

\indent$(4)$ If $G_{e}$ is $C5$, then $|\overline{N[u]}\cap \overline{N[v]}|=5$;

\indent$(5)$ If $G_{e}$ is $C7$ or $C9$, then $2 \leq |\overline{N[u]}\cap \overline{N[v]}|\leq 5$;

\indent$(6)$ If $G_{e}$ is $C8$ or $C11$, then $|\overline{N[u]}\cap \overline{N[v]}|\geq 6$;

\indent$(7)$ If $G_{e}$ is $C10$, then $4 \leq |\overline{N[u]}\cap \overline{N[v]}|\leq 5$;

\indent$(8)$ If $G_{e}$ is $C12$, then $|\overline{N[u]}\cap \overline{N[v]}|\geq 5$;

\indent$(9)$ If $G_{e}$ is $C14$, then $|\overline{N[u]}\cap \overline{N[v]}|\geq 4$.

We show only Claim 2 for configurations $C1$ and $C13$.
The proofs in the other configurations are similar and thus omitted.

If $G_{e}$ is $C1$, then $\overline{N[u]}\supseteq \{v_{1}, v_{2}, v_{3}, v_{4}, v_{5}, v_{6}\}$
and $\overline{N[u]}$ contains at most one vertex in $S_{e}$, which possibly belongs to $\overline{N[v]}$.
Since $C_{v}$ is factor-critical graph, $v$ has at least two neighbors in
$\{v_{1}, v_{2}, v_{3}, v_{4}, v_{5}, v_{6}\}$.
So $|\overline{N[u]}\cap \overline{N[v]}|\leq 5$.

If $G_{e}$ is $C13$, then it is easy to see that $\overline{N[u]}\supseteq \{v_{1}, v_{2}, v_{3}, v_{4}, v_{5}\}$.
Since $vv_{1}, vv_{2}\in$ $E(G)$, $|\overline{N[u]}\cap \overline{N[v]}|\leq 5$.
Moreover, $\{v_{3}, v_{4}, v_{5}\}$$\subseteq \overline{N[u]}\cap \overline{N[v]}$.
So $3 \leq |\overline{N[u]}\cap \overline{N[v]}|\leq 5$. Further, $\{v_{3}, v_{4}, v_{5}\}$ is an
independent set of $G$.

\smallskip
By Claim 1, there are fourteen configurations to discuss.
Next we will complete the entire proof by obtaining a contradiction to each configuration.

\smallskip
{\textbf{Case 1.} $G_{e}$ is $C1$.}

Let $M$ be a perfect matching of $G-S_{e}$. Then $e=uv\in M$. We may assume that
$v_{1}v_{2}, v_{3}v_{4}, v_{5}v_{6}\in M$. We apply Claim 1 to another edge $e'=uu_{1}$ (see $C1$ of Fig. 1).
That is, there exists $S_{e'}\subseteq V(G)-\{u, u_{1}\}$ with $|S_{e'}|=n-10$ such that
$G_{e'}=G-e'-S_{e'}$ is one of Configurations $C1$ to $C14$ relative to edge $e'$. Clearly,
$\overline{N[u]}\cap \overline{N[u_{1}]}=\{v_{1}, v_{2}, v_{3}, v_{4}, v_{5}, v_{6}\}$,
which are paired perfectly under $M$.
By Claim 2, $G_{e'}$ must not be $C1$, $C2$, $C3$, $C4$, $C5$, $C6$, $C7$, $C9$, $C10$ or $C13$.
For the remaining configurations $C8$, $C11$, $C12$ and $C14$,
we cannot find three independent edges in the subgraph induced by the common non-neighborhoods
of $u$ and $u_{1}$ if $G_{e'}$ is $C8$, $C11$, $C12$ or $C14$, a contradiction.

\smallskip
{\textbf{Case 2.} $G_{e}$ is $C4$.}

For a perfect matching $M$ of $G-S_{e}$, also we may assume that
$v_{1}v_{2}, v_{3}v_{4}, v_{5}v_{6}, av_{7}\in M$.
We claim that $av_{5}, av_{6}\in E(G)$. Otherwise, say $av_{5}\notin E(G)$.
Then we consider edge $v_{5}v_{7}$. Clearly,
$\overline{N[v_{5}]}\cap \overline{N[v_{7}]}=\{u, v, v_{1}, v_{2}, v_{3}, v_{4}\}$ and
$uv, v_{1}v_{2}, v_{3}v_{4}$ are three independent edges. By a similar discussion with Case 1,
$G-v_{5}v_{7}-S_{v_{5}v_{7}}$ is not any one of Configurations $C1$ to $C14$ for any
$S_{v_{5}v_{7}}\subseteq V(G)-\{v_{5}, v_{7}\}$ with $|S_{v_{5}v_{7}}|=n-10$, which contradicts Claim 1.

Consider edge $e'=ua$. Obviously,
$\overline{N[u]}\cap \overline{N[a]}\subseteq \{v_{1}, v_{2}, v_{3}, v_{4}\}$.
By Claim 2 and Case 1, $G_{e'}$ is not $C1$, $C3$, $C5$, $C8$, $C11$ or $C12$.
Since $v_{1}v_{2}, v_{3}v_{4}$ are two independent edges,
$G_{e'}$ is not $C10$, $C13$ or $C14$. So $G_{e'}$ is $C2$, $C4$, $C6$, $C7$ or $C9$.

If $G_{e'}$ is $C2$, then $C_{u}$ and $C_{a}$ are two components of $G_{e'}$ with five vertices.
Since $u$ is adjacent to each vertex in $S_{e}\cup \{v\}$, $C_{a}$ contains four vertices among
$\{v_{1}, v_{2}, v_{3}, v_{4}, v_{5}, v_{6}, v_{7}\}$ forming two independent edges
(using the same vertex labeling as $C4$ relative to $e$).
If $C_{a}$ contains exactly two vertices in $\{v_{5}, v_{6}, v_{7}\}$,
say $v_{5}$ and $v_{6}$, then $C_{a}$ contains a pair of adjacent vertices, say $v_{1}$ and $v_{2}$.
Then $v_{7}\in \overline{N[v_{1}]}=\{v_{5}, v_{6}\}\cup V(C_{u})$.
So $v_{7}\in V(C_{u})$ but $v_{6}v_{7}\in E(G)$, a contradiction.
Thus $v_{1}, v_{2}, v_{3}, v_{4}\in V(C_{a})$.
Then $\overline{N[v_{1}]}\supseteq \{v_{5}, v_{6},$ $v_{7}\}\cup V(C_{u})$.
Since $av_{5}, av_{6}, av_{7}\in E(G)$, $v_{5}, v_{6}, v_{7} \notin V(C_{u})$.
So $d_{G}(v_{1})\leq n-9$, a contradiction.

If $G_{e'}$ is $C4$, then we may assume that $G[\{v_{1}, v_{2}, v_{3}\}]$
is a nontrivial odd component of $C4-X$ as $av_{5}, av_{6}, av_{7}\in E(G)$.
It follows that $a$ (resp. $u$) belongs to the trivial
(resp. nontrivial) odd component of $C4-X$ (see Fig. 2).
Otherwise, $v_{4}\in V(C_{a})$ but $v_{3}v_{4}\in E(G)$, a contradiction.
Then $\overline{N[v_{1}]}\supseteq \{a, v_{5}, v_{6}, v_{7}\}\cup V(C_{u})$.
Since $av_{5}, av_{6}, av_{7}\in E(G)$, $v_{5}, v_{6},$ $v_{7} \notin V(C_{u})$.
So $d_{G}(v_{1})\leq n-10$, a contradiction.

\begin{figure}[h]
\centering
\includegraphics[height=2.8cm,width=8cm]{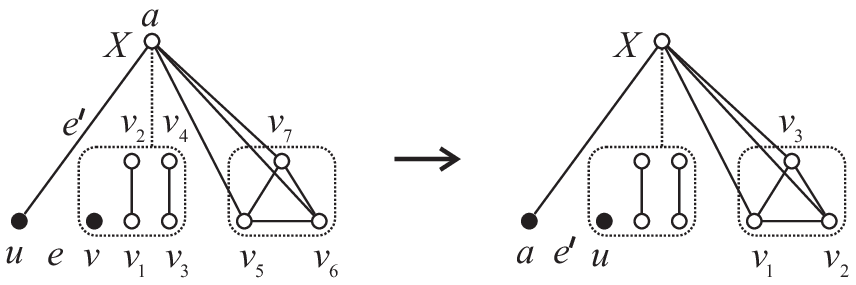}
\caption{\label{tu-1} $G_{e'}$ is $C4$ by applying Claim 1 to edge $e'$.}
\end{figure}

If $G_{e'}$ is $C6$, then three vertices among $\{v_{1}, v_{2}, v_{3}, v_{4}\}$ would induce
a component $K_{3}$ of $C6-X$, say $G[\{v_{1}, v_{2}, v_{3}\}]$.
Thus $C_{a}$ may be $G[\{a, v_{5}, v_{6}\}]$.
Assume that $G[\{u, x_{1}, x_{2}\}]$ is another component of $C6-X$.
Then $\overline{N[v_{5}]}\supseteq \{u, v, v_{1},$ $v_{2}, v_{3}, v_{4}, x_{1}, x_{2}\}$.
Since $ux_{1}, ux_{2}\in E(G)$ and $uv_{4}\notin E(G)$, $v_{4}\notin \{x_{1}, x_{2}\}$.
Because $v$ has at least one neighbor in $\{v_{1}, v_{2}, v_{3}\}$, $v\notin \{x_{1}, x_{2}\}$.
So $d_{G}(v_{5})\leq n-9$, a contradiction.

If $G_{e'}$ is $C7$, then we may choose $v_{1}, v_{3}$ as two trivial odd components of $C7-X$.
It follows that $a$ (resp. $u$) belongs to the trivial (resp. nontrivial) odd component of $C7-X$ (see Fig. 3).
Otherwise, $v_{2}$ or $v_{4}\in V(C_{a})$ but $v_{1}v_{2}, v_{3}v_{4}\in E(G)$, a contradiction.
Then $\{v_{5}, v_{6}, v_{7}\}\subseteq \overline{N[v_{1}]}=\{a, v_{3}\}\cup V(C_{u})$.
So $v_{5}, v_{6}, v_{7}\in V(C_{u})$ but $av_{5}, av_{6}, av_{7}\in E(G)$, a contradiction.

\begin{figure}[h]
\centering
\includegraphics[height=2.8cm,width=8cm]{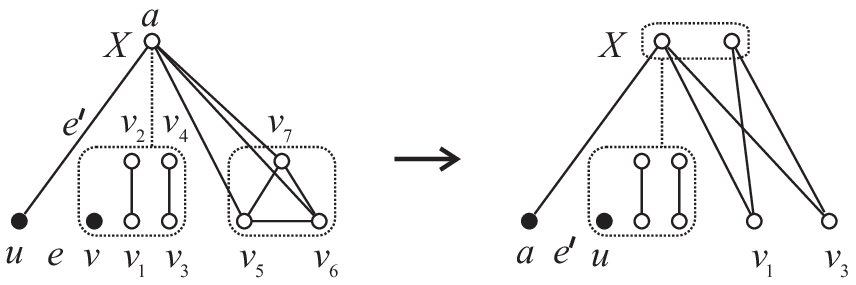}
\caption{\label{tu-5} $G_{e'}$ is $C7$ by applying Claim 1 to edge $e'$.}
\end{figure}

If $G_{e'}$ is $C9$, similarly, we may assume $v_{1}, v_{3}$ as the two trivial odd components of $C9-X$,
and $G[\{a, v_{5}, v_{6}\}]$ and $G[\{u, x_{1}, x_{2}\}]$ as the other two odd components.
Then $\overline{N[v_{1}]}\supseteq \{u, x_{1}, x_{2}, a, v_{3}, v_{5}, v_{6}, v_{7}\}$.
Since $ux_{1}, ux_{2}\in E(G)$ and $uv_{7}\notin E(G)$, $v_{7}\notin \{x_{1}, x_{2}\}$.
So $d_{G}(v_{1})\leq n-9$, a contradiction.

\smallskip
{\textbf{Case 3.} $G_{e}$ is $C9$.}

Take an edge $e'=v_{1}a_{1}$. We apply Claim 1 to $e'$.
Clearly, $\overline{N[v_{1}]}\cap \overline{N[a_{1}]}\subseteq \{u, v, u_{1},$ $u_{2}, v_{3}, v_{4}\}$,
which induces two triangles with an edge between them.
By Claim 2 and Cases 1 and 2, it is obvious that $G_{e'}$ is not $C1$, $C3$ or $C4$.
For $C5$ and $C8$, $G[\overline{N[v_{1}]}\cap \overline{N[a_{1}]}]$ contains
a factor-critical subgraph with 5-vertices.
For $C11$, it consists of two disjoint triangles.
For $C12$, $C13$ and $C14$, it contains an independent set of three vertices.
Such situations would be impossible. So there are five remaining cases to discuss.

If $G_{e'}$ is $C2$, then $C_{a_{1}}$ contains four vertices in $\overline{N[v_{1}]}$ forming
two independent edges. We may assume that $u_{1}, u_{2}\in V(C_{a_{1}})$. Then at least one of
$v, v_{3}$ and $v_{4}$ belongs to $V(C_{a_{1}})$, say $v\in V(C_{a_{1}})$.
Thus $v_{2}\in \overline{N[v]}=\{u_{1}, u_{2}\}\cup V(C_{v_{1}})$.
So $v_{2}\in V(C_{v_{1}})$ but $a_{1}v_{2}\in E(G)$, a contradiction.

If $G_{e'}$ is $C6$, then let $G[\{v, v_{3}, v_{4}\}]$ be a component of $C6-X$
as $G[\overline{N[v_{1}]}\cap \overline{N[a_{1}]}]$ contains a $K_{3}$ in $C6$.
Thus $C_{a_{1}}$ must be $G[\{a_{1}, u_{1}, u_{2}\}]$. Assume that
$C_{v_{1}}$ is $G[\{v_{1},$ $x_{1}, x_{2}\}]$ (see Fig. 4).
Thus $\overline{N[v_{3}]}\supseteq \{a_{1}, v_{1}, u_{1}, u_{2}, x_{1}, x_{2},$ $u, v_{2}\}$.
Since $v_{1}x_{1}$, $v_{1}x_{2}\in E(G)$ and $v_{1}v_{2}$, $v_{1}u\notin E(G)$,
$x_{1}, x_{2}\notin \{u, v_{2}\}$. So $d_{G}(v_{3})$$\leq n-9$, a contradiction.

\begin{figure}[h]
\centering
\includegraphics[height=2.8cm,width=8cm]{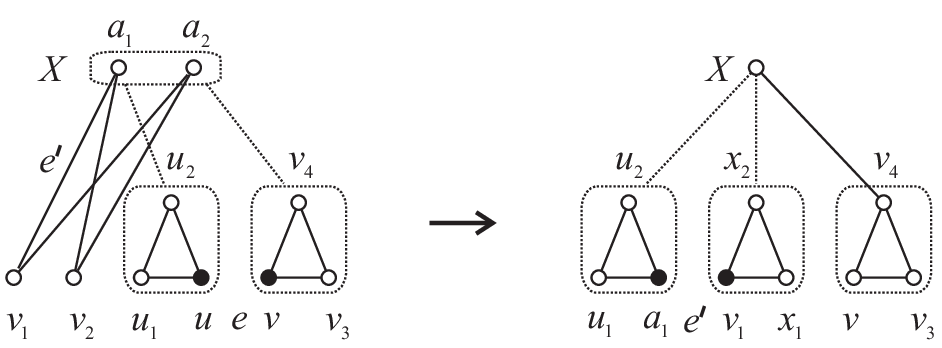}
\caption{\label{tu-6} $G_{e'}$ is $C6$ by applying Claim 1 to edge $e'$.}
\end{figure}

If $G_{e'}$ is $C7$, then we may assume $u, v_{3}$ as two trivial odd components of $C7-X$.
It follows that $a_{1}$ is the third trivial component of $C7-X$.
Otherwise, $v_{1}$ is the third trivial component, and $u_{1}$ or $u_{2}\in V(C_{a_{1}})$
but $uu_{1}, uu_{2}\in E(G)$, a contradiction.
Then $v_{2}\in \overline{N[v_{3}]}=\{a_{1}, u\}\cup V(C_{v_{1}})$.
So $v_{2}\in V(C_{v_{1}})$ but $a_{1}v_{2}\in E(G)$, a contradiction.

If $G_{e'}$ is $C9$, similarly we may choose $u, v_{3}$ as the two trivial components of $C9-X$.
Thus $C_{a_{1}}$ is either $G[\{a_{1}, u_{1}, u_{2}\}]$ or $G[\{a_{1}, v, v_{4}\}]$.
But $uu_{1}, uu_{2}, vv_{3}, v_{3}v_{4}\in E(G)$, a contradiction.

If $G_{e'}$ is $C10$, then $G[\overline{N[v_{1}]}\cap \overline{N[a_{1}]}]$ contains a $K_{1}$
and a $K_{3}$ in $C10$, which are disjoint. So we may assume $u_{1}$ as a trivial
component and $G[\{v, v_{3}, v_{4}\}]$ is a nontrivial component of $C10-X$.
It follows that $a_{1}$ (resp. $v_{1}$) belongs to the trivial (resp. nontrivial) component of $C10-X$.
Otherwise, $u$ or $u_{2}\in V(C_{a_{1}})$ but $uu_{1}, u_{1}u_{2}\in E(G)$, a contradiction.
Assume that $C_{v_{1}}$ is $G[\{v_{1}, x_{1}, x_{2}\}]$.
Then $\overline{N[u_{1}]}\supseteq\{a_{1}, x_{1}, x_{2}, v, v_{1}, v_{2}, v_{3}, v_{4}\}$.
Since $v_{1}x_{1}$, $v_{1}x_{2}\in E(G)$ and $v_{1}v_{2}\notin E(G)$,
$v_{2}\notin \{x_{1}, x_{2}\}$. So $d_{G}(u_{1})\leq n-9$, a contradiction.

\smallskip
{\textbf{Case 4.} $G_{e}$ is $C5$.}

For a perfect matching $M$ of $G-S_{e}$, we may assume that
$av_{3}, v_{4}v_{5}, v_{6}v_{7}\in M$. We discuss the three configurations of $C5$
as shown in Fig. 5.
(By symmetry, $v_{1}$ and $v_{2}$ are equivalent, and the dotted edge is an optional edge.)

\begin{figure}[h]
\centering
\includegraphics[height=3.3cm,width=11cm]{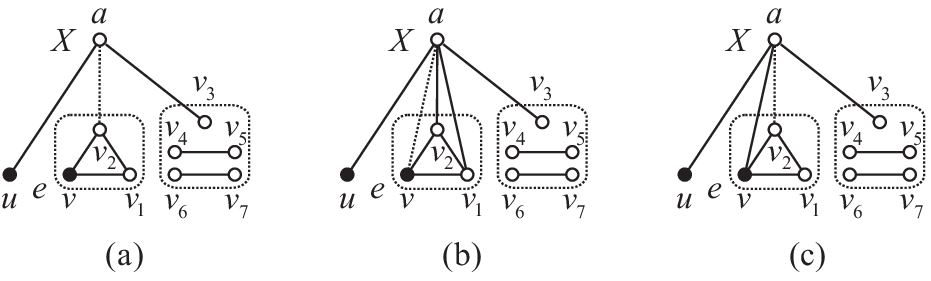}
\caption{\label{tu-02} The three configurations of $C5$.}
\end{figure}

{\textbf{Subcase 4.1.} $av, av_{1}\notin E(G)$ and $av_{2}$ is an optional edge (see Fig. 5 (a)).}

Consider edge $vv_{1}$. Then $\overline{N[v]}\cap \overline{N[v_{1}]}=\{a, v_{3}, v_{4}, v_{5}, v_{6}, v_{7}\}$
and $av_{3}, v_{4}v_{5}, v_{6}v_{7}$ are three independent edges. We apply Claim 1 to $vv_{1}$.
Then the proof is similar to Case 1.

{\textbf{Subcase 4.2.} $av_{1}, av_{2}\in E(G)$ and $av$ is an optional edge (see Fig. 5 (b)).}

Let $e'=ua$. Then $\overline{N[u]}\cap \overline{N[a]}\subseteq\{v_{4}, v_{5}, v_{6}, v_{7}\}$
and $v_{4}v_{5}, v_{6}v_{7}$ are two independent edges.
By Claim 2 and Cases 1 to 3, it is obvious that $G_{e'}$ is not $C1$, $C3$, $C4$, $C5$,
$C8$, $C9$, $C11$ or $C12$.
For $C10$, there are not two independent edges in the subgraph induced by the common
non-neighborhoods of $u$ and $a$.
For $C13$ and $C14$, it contains at least three independent vertices.
Both of them contradict that $v_{4}v_{5}, v_{6}v_{7}$ are two independent edges.
Then there are three remaining cases to discuss.

If $G_{e'}$ is $C2$, then $C_{a}$ contains four vertices among
$\{v_{1}, v_{2}, v_{3}, v_{4}, v_{5}, v_{6}, v_{7}\}$ forming two independent edges.
We need to consider the two subcases depending on whether $v_{1}, v_{2}\in V(C_{a})$ or not.

When $v_{1}, v_{2}\in V(C_{a})$, we may assume that $v_{3}, v_{4}\in V(C_{a})$.
Then $\{v_{5}, v_{6}, v_{7}\}\subseteq \overline{N[v_{1}]}=\{v_{3}, v_{4}\}\cup V(C_{u})$.
So $v_{5}, v_{6}, v_{7}\in V(C_{u})$ but $v_{4}v_{5}\in E(G)$, a contradiction.

When $v_{1}, v_{2}\notin V(C_{a})$, we assume that $v_{3}, v_{4}, v_{5}, v_{6}\in V(C_{a})$
as show in Fig. 6. Then $\overline{N[v_{3}]}\supseteq \{v, v_{1}, v_{2}\}\cup V(C_{u})$.
Since $av_{1}, av_{2}\in E(G)$, $v_{1}, v_{2}\notin V(C_{u})$
and $v\in V(C_{u})$. So $av\notin E(G)$. Assume that $C_{u}=G[\{u, v, u_{1}, u_{2}, u_{3}\}]$
and $vu_{1}, u_{2}u_{3}$ are two independent edges.
Thus $\overline{N[v_{i}]}=\{v_{1}, v_{2}\}\cup$ $V(C_{u})$ for $i=3, 4, 5, 6$,
which implies that $C_{a}$ is a $K_{5}$.

\begin{figure}[h]
\centering
\includegraphics[height=2.5cm,width=8cm]{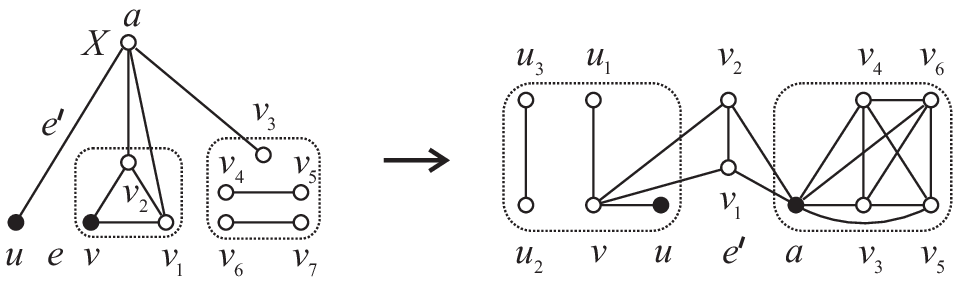}
\caption{\label{tu-2} $G_{e'}$ is $C2$ by applying Claim 1 to $e'$.}
\end{figure}

Now consider edge $av_{3}$. Since $au, av_{1}, av_{2}\in E(G)$,
$\overline{N[a]}\cap \overline{N[v_{3}]}=\{v, u_{1}, u_{2}, u_{3}\}$.
By Claim 2 and Cases 1 to 3, we obtain that $G-av_{3}-S_{av_{3}}$ is not $C1$, $C3$, $C4$,
$C5$, $C8$, $C9$, $C10$, $C11$, $C12$, $C13$ or $C14$ again.

If $G-av_{3}-S_{av_{3}}$ is $C2$, then $C_{v_{3}}$ contains at least one vertex in $\{v, u_{1}, u_{2}, u_{3}\}$
as $|\overline{N[a]}|\leq 7$ and $|\overline{N[a]}\cap \overline{N[v_{3}]}|=4$, say $u_{1}\in V(C_{v_{3}})$.
Then $\overline{N[u_{1}]}\supseteq\{v_{3}, v_{4}, v_{5}, v_{6}\}\cup V(C_{a})$.
Since $v_{3}v_{4}, v_{3}v_{5}, v_{3}v_{6}\in E(G)$, $v_{4}, v_{5}, v_{6} \notin V(C_{a})$.
So $d_{G}(u_{1})\leq n-10$, a contradiction.

If $G-av_{3}-S_{av_{3}}$ is $C6$, then $C_{a}$ must be $G[\{a, v_{1}, v_{2}\}]$ as $uv_{1}, uv_{2}\notin E(G)$.
Since $vv_{1}, vv_{2}\in E(G)$, $G[\{u_{1}, u_{2}, u_{3}\}]$ is another component of $C6-X$.
Hence $\overline{N[v_{1}]}\supseteq$ $\{u, v_{3},$ $v_{4}, v_{5}, v_{6}, u_{1},$ $u_{2}, u_{3}\}$.
So $d_{G}(v_{1})\leq n-9$, a contradiction.

If $G-av_{3}-S_{av_{3}}$ is $C7$, then we may choose $u_{1}, u_{3}$ as
two trivial odd components of $C7-X$. It follows that $a$ (resp. $v_{3}$) belongs to the trivial
(resp. nontrivial) odd component of $C7-X$.
Otherwise, $v$ or $u_{2}\in V(C_{a})$ but $vu_{1}, u_{2}u_{3}\in E(G)$, a contradiction.
Then $\{v_{4}, v_{5}, v_{6}\}\subseteq \overline{N[u_{1}]}=\{a, u_{3}\}\cup V(C_{v_{3}})$.
So $v_{4}, v_{5}, v_{6}\in V(C_{v_{3}})$ but $av_{4}, av_{5}, av_{6}\in E(G)$, a contradiction.

The above discussions imply that $G_{e'}$ is not $C2$.

\smallskip

If $G_{e'}$ is $C6$, then we assume that $G[\{v_{4}, v_{5}, v_{6}\}]$ and
$G[\{u, x_{1}, x_{2}\}]$ are two odd components of $C6-X$. Thus $C_{a}$
must be $G[\{a, v_{1}, v_{2}\}]$.
Hence $\overline{N[v_{1}]}\supseteq \{u, v_{3}, v_{4}, v_{5}, v_{6}, v_{7},$ $x_{1}, x_{2}\}$.
Since $ux_{1}, ux_{2}\in E(G)$ and $uv_{3}, uv_{7}\notin E(G)$,
$x_{1}, x_{2}\notin \{v_{3}, v_{7}\}$. So $d_{G}(v_{1})\leq n-9$, a contradiction.

If $G_{e'}$ is $C7$, then we may choose $v_{4}, v_{6}$ as two
trivial components $C7-X$. It follows that $a$ is the third trivial component of $C7-X$.
Otherwise, $u$ is the third trivial component, and $v_{5}$ or $v_{7}\in V(C_{a})$
but $v_{4}v_{5}, v_{6}v_{7}\in E(G)$, a contradiction.
Then $\{v_{1}, v_{2}\}$ $\subseteq \overline{N[v_{4}]}=\{a, v_{6}\}\cup V(C_{u})$.
So $v_{1}, v_{2}\in V(C_{u})$ but $av_{1}, av_{2}\in E(G)$, a contradiction.

{\textbf{Subcase 4.3.} $av\in E(G), av_{1}\notin E(G)$ and $av_{2}$ is an optional edge (see Fig. 5 (c)).}

Take an edge $e'=vv_{1}$. It is easy to see that
$\overline{N[v]}\cap \overline{N[v_{1}]}=\{v_{3}, v_{4}, v_{5}, v_{6}, v_{7}\}$.
By Claim 2 and Cases 1 to 3, $G_{e'}$ is not $C1$, $C3$, $C4$, $C8$, $C9$ or $C11$.
Since $G[\overline{N[v]}\cap \overline{N[v_{1}]}]$ is factor-critical,
it does contain a $K_{1}$ and a $K_{3}$ which are disjoint as induced subgraphs.
So $G_{e'}$ is not $C10$. For $C12$, $C13$ and $C14$,
it also does not contain an independent set with three vertices.
So there are four remaining cases to discuss.

If $G_{e'}$ is $C2$, then $u, a\in V(C_{v})$. Assume that $v_{3}, v_{4}\in V(C_{v})$.
Then $\overline{N[v_{3}]}\supseteq \{u, v, v_{2}\} \cup V(C_{v_{1}})$.
Since $vv_{2}\in E(G)$, $v_{2}\notin V(C_{v_{1}})$. So $d_{G}(v_{3})\leq n-9$, a contradiction.

If $G_{e'}$ is $C6$, then $C_{v}$ must be $G[\{a, u, v\}]$ as $v$ has only
two neighbors $u$ and $a$ in $\overline{N[v_{1}]}$. Assume that the other two
components of $C6-X$ are $G[\{v_{4}, v_{5}, v_{6}\}]$ and $G[\{v_{1}, x_{1}, x_{2}\}]$.
Hence $\{x_{1}, x_{2}\}\subseteq \overline{N[u]}= \{v_{1}, v_{2}, v_{3}, v_{4}, v_{5},$
$v_{6}, v_{7}\}$. So $x_{1}, x_{2}\in \{v_{2}, v_{3}, v_{7}\}$
but $v_{1}x_{1}, v_{1}x_{2}\in E(G)$ and $v_{1}v_{3}, v_{1}v_{7}\notin E(G)$, a contradiction.

If $G_{e'}$ is $C7$, then we may assume that $v_{4}, v_{6}$
as two trivial odd components of $C7-X$. It follows that $v$ (resp. $v_{1}$)
belongs to the trivial (resp. nontrivial) odd component of $C7-X$.
Otherwise, $v_{5}$ or $v_{7}\in V(C_{v})$ but $v_{4}v_{5}, v_{6}v_{7}\in E(G)$, a contradiction.
Then $\{u, v_{2}\}$ $\subseteq \overline{N[v_{4}]}=\{v, v_{6}\}\cup V(C_{v_{1}})$.
So $u, v_{2}\in V(C_{v_{1}})$ but $uv, vv_{2}\in E(G)$, a contradiction.

If $G_{e'}$ is $C5$, then it still has configuration as show in Fig. 5 (c) by Subcases 4.1 and 4.2.
So $G[\{v_{3}, v_{4}, v_{5}, v_{6}, v_{7}\}]$ is an odd component of $C5-X$.
Then $v$ (resp. $v_{1}$) belongs to the trivial (resp. nontrivial)
odd component. Otherwise, $C_{v}$ would be $G[\{u, v, a\}]$
but $av_{3}\in E(G)$, a contradiction. Assume that $C_{v_{1}}$ is $G[\{v_{1}, x_{1}, x_{2}\}]$
(see Fig. 7).

\begin{figure}[h]
\centering
\includegraphics[height=2.8cm,width=12.5cm]{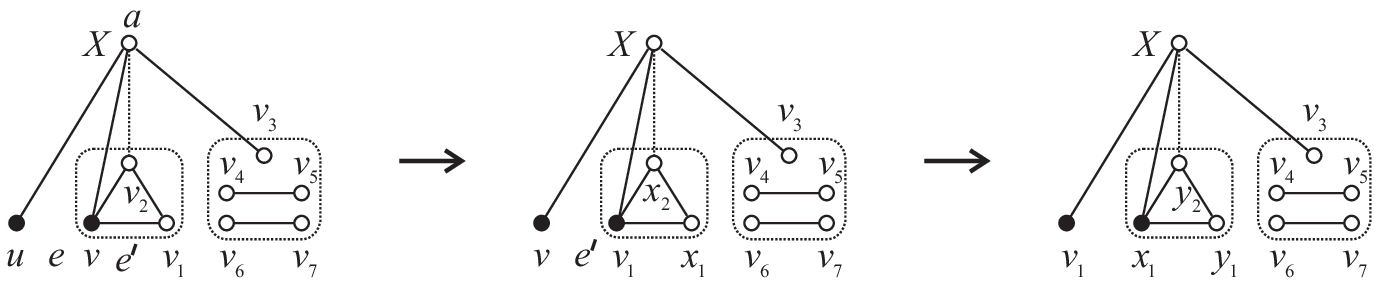}
\caption{\label{tu-3} $G_{e'}$ and $G_{v_{1}x_{1}}$ are still $C5$ by applying Claim 1 to $e'$ and
$v_{1}x_{1}$, respectively.}
\end{figure}

Now consider edge $v_{1}x_{1}$ and apply Claim 1 to $v_{1}x_{1}$.
Then $\overline{N[v_{1}]}\cap \overline{N[x_{1}]}=\{v_{3}, v_{4}, v_{5},$ $v_{6}, v_{7}\}$.
By the former discussions in Subcase 4.3, $G-v_{1}x_{1}-S_{v_{1}x_{1}}$ would also be $C5$.
Then $G[\{v_{3}, v_{4}, v_{5}, v_{6}, v_{7}\}]$ is still the odd component of $C5-X$.
Similarly, $v_{1}$ (resp. $x_{1}$) belongs to the trivial (resp. nontrivial) odd component of $C5-X$.
Assume that $C_{x_{1}}$ is $G[\{x_{1}, y_{1}, y_{2}\}]$.
Thus $\overline{N[v_{1}]}\supseteq \{u, a, v_{3}, v_{4}, v_{5}, v_{6}, v_{7}, y_{1}, y_{2}\}$.
Since $av_{3}\in E(G)$ and $y_{1}v_{3}, y_{2}v_{3}\notin E(G)$,
$y_{1}, y_{2}\neq a$ (possibly, $y_{1}$ or $y_{2}=u$). So $d_{G}(v_{1})\leq n-9$, a contradiction.

\smallskip
{\textbf{Case 5.} $G_{e}$ is $C2$.}

Let $e'=u_{3}u_{4}$. We apply Claim 1 to $e'$.
Clearly, $|\overline{N[u_{3}]}\cap \overline{N[u_{4}]}|\geq 5$.
We divide the proof into the following three subcases according to
$|\overline{N[u_{3}]}\cap \overline{N[u_{4}]}|$.

{\textbf{Subcase 5.1.} $|\overline{N[u_{3}]}\cap \overline{N[u_{4}]}|=5$.}

Then $\overline{N[u_{3}]}\cap \overline{N[u_{4}]}=\{v, v_{1}, v_{2}, v_{3}, v_{4}\}$.
By Claim 2 and Cases 1 to 4, $G_{e'}$ is not $C1$, $C3$, $C4$, $C5$, $C8$, $C9$, $C10$, $C11$, $C12$,
$C13$ or $C14$.

If $G_{e'}$ is $C2$, then $C_{u_{3}}$ and $C_{u_{4}}$ contain respectively two vertices in
$\{v, v_{1}, v_{2}, v_{3},$ $v_{4}\}$, say $v_{1}, v_{2}\in V(C_{u_{3}})$ and $v_{3}, v_{4}\in V(C_{u_{4}})$.
Then $\overline{N[v_{1}]}=\{u, u_{1}, u_{2}, u_{3}, u_{4}, v_{3}, v_{4}\}$.
So two of $u, u_{1}$ and  $u_{2}$ belong to $V(C_{u_{4}})$
but no edges join $\{v_{3}, v_{4}\}$ and $\{u, u_{1}, u_{2}, u_{4}\}$,
contradicting that $C_{u_{4}}$ is connected.

If $G_{e'}$ is $C6$, then we may assume that $G[\{v_{1}, v_{2}, v_{3}\}]$,
$G[\{u_{3}, x_{1}, x_{2}\}]$ and $G[\{u_{4}, x_{3},$ $x_{4}\}]$ are three components of $C6-X$ as show in Fig. 8.
Hence $\overline{N[v_{1}]}\supseteq\{u, u_{1}, u_{2}, u_{3}, u_{4}, x_{1}, x_{2},$ $x_{3}, x_{4}\}$.
Since $d_{G}(v_{1})\geq n-8$, $|\{u, u_{1}, u_{2}\}\cap \{x_{1}, x_{2}, x_{3}, x_{4}\}|\geq 2$,
say $\{u_{1}, u_{2}\}=\{x_{1}, x_{2}\}$ (similarly, $\{u, u_{1}\}=\{x_{1}, x_{3}\}$).
Then $\overline{N[u_{1}]}\supseteq \{v, v_{1}, v_{2}, v_{3}, v_{4}, u_{4}, x_{3}, x_{4}\}$.
Since $u_{4}x_{3}, u_{4}x_{4}$ $\in E(G)$ and $u_{4}v, u_{4}v_{4}\notin E(G)$,
$x_{3}, x_{4}\notin \{v, v_{4}\}$. So $d_{G}(u_{1})\leq n-9$, a contradiction.

\begin{figure}[h]
\centering
\includegraphics[height=2.6cm,width=9cm]{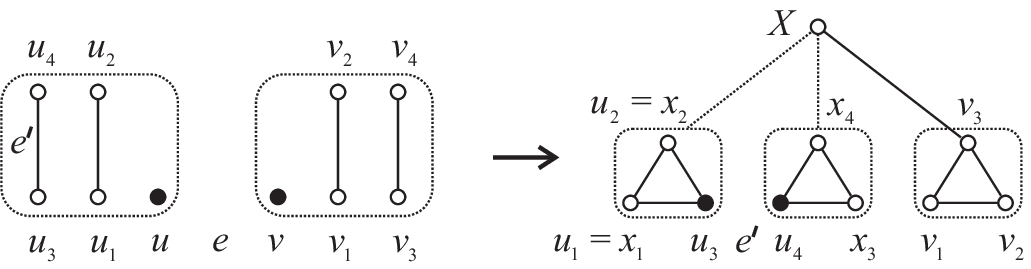}
\caption{\label{tu-4} $G_{e'}$ is $C6$ by applying Claim 1 to $e'$.}
\end{figure}

If $G_{e'}$ is $C7$, without loss of generality, we assume that $u_{3}$ (resp. $u_{4}$)
belong to the trivial (resp. nontrivial) odd component of $C7-X$ and
choose $v, v_{1}$ (similarly, $v_{1}, v_{3}$) as the other two trivial odd components.
Then $\{u, u_{1}, u_{2}\}\subseteq \overline{N[v_{1}]}=\{v, u_{3}\}\cup V(C_{u_{4}})$.
So $u, u_{1}, u_{2}\in V(C_{u_{4}})$,
contradicting that $G[\{u, u_{1}, u_{2}, u_{3}, u_{4}\}]$ is factor-critical.

{\textbf{Subcase 5.2.} $|\overline{N[u_{3}]}\cap \overline{N[u_{4}]}|=6$.}

We consider the four situations of $|\overline{N[u_{3}]}\cap \overline{N[u_{4}]}|=6$ as show in Fig. 9.
(The dotted edge is an optional edge and black vertices are the vertices
in $\overline{N[u_{3}]}\cap \overline{N[u_{4}]}$.)

\begin{figure}[h]
\centering
\includegraphics[height=2.6cm,width=16cm]{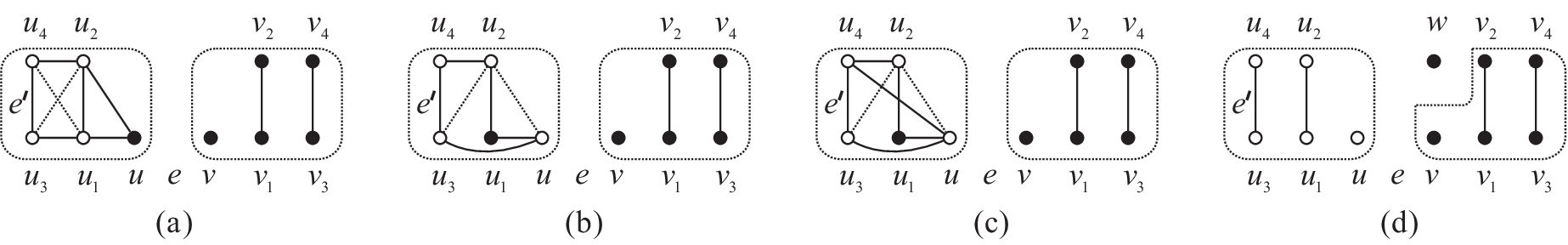}
\caption{\label{tu-03} The four situations of $|\overline{N[u_{3}]}\cap \overline{N[u_{4}]}|=6$.}
\end{figure}

(a) Then $\overline{N[u_{3}]}\cap \overline{N[u_{4}]}=\{u, v, v_{1}, v_{2}, v_{3}, v_{4}\}$
and $uv, v_{1}v_{2}, v_{3}v_{4}$ are three independent edges.
We apply Claim 1 to edge $u_{3}u_{4}$ and then the proof is similar to Case 1.

(b) If $uu_{2}\notin E(G)$,
then $\overline{N[u_{2}]}\cap \overline{N[u_{4}]}=\{u, v, v_{1}, v_{2}, v_{3}, v_{4}\}$
and $uv, v_{1}v_{2}, v_{3}v_{4}$ are three independent edges.
We consider edge $u_{2}u_{4}$. The proof is similar to the Case 1.

If $uu_{2}\in E(G)$, we need to consider two subcases depending on whether $u_{2}u_{3}\in E(G)$ or not.
If $u_{2}u_{3}\in E(G)$,
then $\overline{N[u_{1}]}\cap \overline{N[u_{2}]}=\{v, v_{1}, v_{2}, v_{3}, v_{4}\}$.
We can give a similar proof as Subcase 5.1 by applying Claim 1 to edge $u_{1}u_{2}$.
Otherwise, we consider edge $uu_{3}$.
Clearly, $\overline{N[u]}\cap \overline{N[u_{3}]}=\{v_{1}, v_{2}, v_{3}, v_{4}\}$.
Then the proof is similar to Subcase 4.2.

(c) For edge $uu_{1}$, we have $\overline{N[u]}\cap \overline{N[u_{1}]}=\{v_{1}, v_{2}, v_{3}, v_{4}\}$.
By applying Claim 1 to $uu_{1}$, we can give a similar proof as Subcase 4.2.

(d) Then $\overline{N[u_{3}]}\cap \overline{N[u_{4}]}=\{v, v_{1}, v_{2}, v_{3}, v_{4}, w\}$,
where $w\in S_{e}$. By Claim 2 and $G[\{v, v_{1},$ $v_{2}, v_{3}, v_{4}\}]$
is factor-critical, $G-e'-S_{e'}$ would only be $C8$.
Thus $w$ is a trivial component and
$G[\{v, v_{1}, v_{2}, v_{3}, v_{4}\}]$ is the nontrivial component of $C8-X$.
Then $\overline{N[w]}=\{v, v_{1}, v_{2}, v_{3},$ $v_{4}, u_{3}, u_{4}\}$.
So $uw, u_{1}w, u_{2}w\in E(G)$.
Since $|\overline{N[u_{3}]}|\leq 7$ and $|\overline{N[u_{4}]}|\leq 7$,
$u_{3}$ and $u_{4}$ have at least two neighbors in $\{u, u_{1}, u_{2}\}$.
That is, $u_{3}$ and $u_{4}$ have at least one common neighbor in $\{u, u_{1}, u_{2}\}$.

If $uu_{3}, uu_{4}\in E(G)$, we consider edge $uw$
and $\overline{N[u]}\cap \overline{N[w]}=\{v_{1}, v_{2}, v_{3}, v_{4}\}$.
The proof is similar to Subcase 4.2. Otherwise, say $u_{1}u_{3}, u_{1}u_{4}\in E(G)$.
Then we consider edge $u_{1}w$ and
$\overline{N[u_{1}]}\cap \overline{N[w]}=\{v, v_{1}, v_{2}, v_{3}, v_{4}\}$.
The proof is analogous to the corresponding Subcase 5.1.

{\textbf{Subcase 5.3.} $|\overline{N[u_{3}]}\cap \overline{N[u_{4}]}|=7$.}

We consider the five situations of $|\overline{N[u_{3}]}\cap \overline{N[u_{4}]}|=7$ as show in Fig. 10.
(The black vertices are the vertices of $\overline{N[u_{3}]}\cap \overline{N[u_{4}]}$,
where $w, w_{1}, w_{2}\in S_{e}$.)

\begin{figure}[h]
\centering
\includegraphics[height=5.6cm,width=13cm]{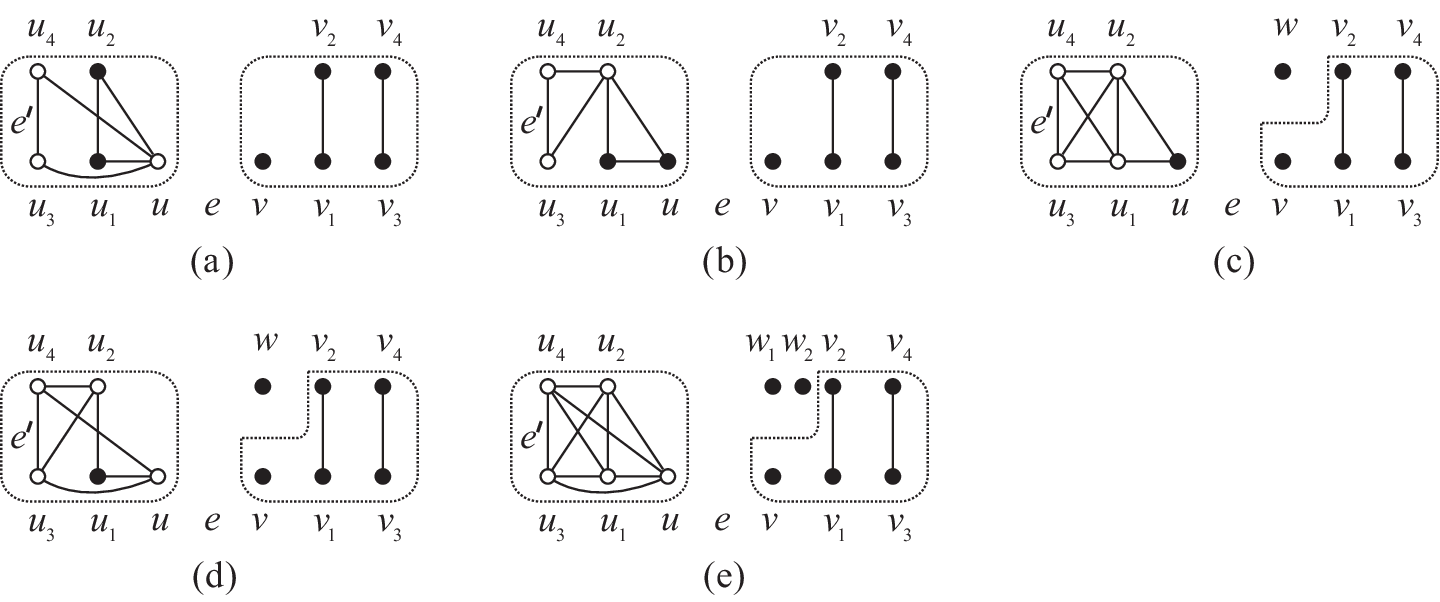}
\caption{\label{tu-04} The five situations of $|\overline{N[u_{3}]}\cap \overline{N[u_{4}]}|=7$.}
\end{figure}

(a) or (d) Consider edge $uu_{1}$.
Clearly, $\overline{N[u]}\cap \overline{N[u_{1}]}=\{v_{1}, v_{2}, v_{3}, v_{4}\}$
and $v_{1}v_{2}, v_{3}v_{4}$ are two independent edges.
The proof is similar to Subcase 4.2 by using Claim 1 to $uu_{1}$.

(b) For edge $uu_{1}$,
$\overline{N[u]}\cap \overline{N[u_{1}]}=\{u_{3}, u_{4}, v_{1}, v_{2}, v_{3}, v_{4}\}$
and $u_{3}u_{4}, v_{1}v_{2}, v_{3}v_{4}$ are three independent edges.
The proof is similar to Case 1 by using Claim 1 to $uu_{1}$.

(c) Take an edge $u_{1}u_{3}$. We apply Claim 1 to $u_{1}u_{3}$.

If $u_{1}w\in E(G)$, then $\overline{N[u_{1}]}\cap \overline{N[u_{3}]}=\{v, v_{1}, v_{2}, v_{3}, v_{4}\}$.
The proof is similar to Subcase 5.1.
Otherwise, $\overline{N[u_{1}]}\cap \overline{N[u_{3}]}=\{w, v, v_{1}, v_{2}, v_{3}, v_{4}\}$.
By Claim 2 and $G[\{v, v_{1}, v_{2}, v_{3}, v_{4}\}]$ is factor-critical,
$G-u_{1}u_{3}-S_{u_{1}u_{3}}$ would only be $C8$.
Thus $w$ belongs to a trivial odd component and
$G[\{v, v_{1}, v_{2}, v_{3}, v_{4}\}]$ is the nontrivial odd component of $C8-X$.
Then $\overline{N[w]}\supseteq\{u_{1}, u_{3}, u_{4}, v, v_{1}, v_{2}, v_{3}, v_{4}\}$.
So $d_{G}(w)\leq n-9$, a contradiction.

(e) If $u_{1}w_{1}, u_{1}w_{2}\in E(G)$,
then $\overline{N[u_{1}]}\cap \overline{N[u_{3}]}=\{v, v_{1}, v_{2}, v_{3}, v_{4}\}$.
We apply Claim 1 to edge $u_{1}u_{3}$ and the proof is similar to Subcase 5.1.

If $u_{1}w_{1}\notin E(G)$, $u_{1}w_{2}\in E(G)$,
then $\overline{N[u_{1}]}\cap \overline{N[u_{3}]}=\{w_{1}, v, v_{1}, v_{2}, v_{3}, v_{4}\}$.
We consider edge $u_{1}u_{3}$ and the proof can be given with a similar argument as Subcase 5.3 (c).

Thus $u_{1}w_{1}, u_{1}w_{2}\notin E(G)$. Similarly, $u_{2}w_{1}, u_{2}w_{2}\notin E(G)$.
So $C_{u}$ is a $K_{5}$.

Next we consider edge $uu_{3}$ and apply Claim 1 to $uu_{3}$.

If $uw_{1}, uw_{2}\in E(G)$, then $\overline{N[u]}\cap \overline{N[u_{3}]}=\{v_{1}, v_{2}, v_{3}, v_{4}\}$
and $v_{1}v_{2}, v_{3}v_{4}$ are two independent edges.
The proof is analogous to the corresponding Subcase 4.2.

If $uw_{1}, uw_{2}\notin E(G)$, then $\overline{N[u]}\cap \overline{N[u_{3}]}=$
$\{w_{1}, w_{2}, v_{1}, v_{2}, v_{3}, v_{4}\}$.
By Claim 2, $G-uu_{3}-S_{uu_{3}}$ would be $C8$, $C11$, $C12$ or $C14$.
Since $\overline{N[w_{1}]} \cap \overline{N[w_{2}]}$$\supseteq \{u, u_{1}, u_{2}, u_{3}, u_{4}\}$,
both $\overline{N[w_{1}]}$ and $\overline{N[w_{2}]}$ contain at most two vertices
in $V(G)\setminus V(C_{u})$.
It is easy to verify that $G-uu_{3}-S_{uu_{3}}$ can not be $C8$, $C11$, $C12$ or $C14$,
which contradicts Claim 1.

Thus we may assume that $uw_{1}\notin E(G)$ and $uw_{2}\in E(G)$.
So $\overline{N[u]}\cap \overline{N[u_{3}]}=\{w_{1}, v_{1}, v_{2}, v_{3}, v_{4}\}$
and $v_{1}v_{2}, v_{3}v_{4}$ are two independent edges.
Since $\overline{N[w_{1}]}$$\supseteq \{u, u_{1}, u_{2}, u_{3},$ $u_{4}\}$,
$w_{1}$ has at least two neighbors in $\{v_{1}, v_{2}, v_{3}, v_{4}\}$.
By Claim 2 and Cases 1 to 4, $G-uu_{3}-S_{uu_{3}}$ is not $C1$, $C3$, $C4$, $C5$, $C8$, $C9$,
$C11$, $C12$ or $C14$. So there are five remaining cases to discuss.

If $G-uu_{3}-S_{uu_{3}}$ is $C2$, then we have $v, w_{2}\in V(C_{u})$ as $uv, uw_{2}\in E(G)$.
Thus $\overline{N[w_{2}]}\supseteq\{u_{1}, u_{2}, u_{4}\}\cup V(C_{u_{3}})$.
Since $uu_{1}, uu_{2}, uu_{4}\in E(G)$, $u_{1}, u_{2}, u_{4}\notin$ $V(C_{u_{3}})$.
So $d_{G}(w_{2})\leq n-9$, a contradiction.

If $G-uu_{3}-S_{uu_{3}}$ is $C6$, then $C_{u}$ must be $G[\{u, v, w_{2}\}]$.
Assume that $G[\{u_{3}, x_{1}, x_{2}\}]$ is an component and $v_{1}$ belongs to another
component of $C6-X$. Then $\overline{N[v_{1}]}\supseteq\{u_{3}, x_{1}, x_{2},$
$u, u_{1}, u_{2}, u_{4}, w_{2}\}$. Since $uu_{1}, uu_{2}, uu_{4}\in E(G)$ and
$ux_{1}, ux_{2}\notin E(G)$, $x_{1}, x_{2}\notin \{u_{1}, u_{2}, u_{4}\}$.
So $d_{G}(v_{1})\leq n-9$, a contradiction.

If $G-uu_{3}-S_{uu_{3}}$ is $C7$, then we choose $v_{1}, v_{3}$
(similarly, $v_{1}, w_{1}$) as two trivial components of $C7-X$.
It follows that $u$ is the third trivial component of $C7-X$.
Otherwise, $u_{3}$ is the third trivial component and $v_{2}$ or $v_{4}\in V(C_{u})$
but $v_{1}v_{2}, v_{3}v_{4}\in E(G)$, a contradiction.
Then $\{u_{1}, u_{2}, u_{4}\}\subseteq \overline{N[v_{1}]}=\{u, v_{3}\}\cup V(C_{u_{3}})$.
So $u_{1}, u_{2}, u_{4}\in V(C_{u_{3}})$ but $uu_{1}, uu_{2}, uu_{4}\in E(G)$, a contradiction.

If $G-uu_{3}-S_{uu_{3}}$ is $C10$, then $w_{1}$ must belong to a nontrivial component of $C10-X$,
say $G[\{w_{1}, v_{1}, v_{2}\}]$. Otherwise, $w_{1}$ belongs to a trivial component
and then $G[\{v_{1}, v_{2}, v_{3}\}]$ would be a nontrivial component of $C10-X$.
Then $\overline{N[w_{1}]}\supseteq\{u, u_{1}, u_{2}, u_{3}, u_{4}, v_{1}, v_{2}, v_{3}\}$,
a contradiction. Assume that $v_{3}$ is a trivial component of $C10-X$.
Then $\overline{N[v_{3}]}\supseteq \{u, u_{1}, u_{2},$ $u_{3}, u_{4}, v_{1}, v_{2}, w_{1}\}$.
So $d_{G}(v_{3})\leq n-9$, a contradiction.

If $G-uu_{3}-S_{uu_{3}}$ is $C13$, then we assume that $\{v_{1}, v_{3}, w_{1}\}$
is an independent set of $G$. It follows that $u$ (resp. $u_{3}$) belongs to a trivial (resp. nontrivial)
component of $C13-X$. Otherwise, $C_{u}$ is $G[\{u, v, w_{2}\}]$ and
then $\{v, v_{1}, v_{3}\}$ is an independent set of $G$, contradicting that
$G[\{v, v_{1}, v_{2}, v_{3}, v_{4}\}]$ is factor-critical.
Assume that $C_{u_{3}}$ is $G[\{u_{3}, x_{1}, x_{2}\}]$ (see Fig. 11).
Thus $\overline{N[w_{1}]}\supseteq\{u, u_{1}, u_{2}, u_{3}, u_{4}, v_{1}, v_{3}, x_{1}, x_{2}\}$.
Since $uu_{1}, uu_{2}, uu_{4}\in E(G)$ and $ux_{1}, ux_{2}\notin E(G)$,
$x_{1}, x_{2}\notin \{u_{1}, u_{2}, u_{4}\}$. So $d_{G}(w_{1})\leq n-10$, a contradiction.

\begin{figure}[h]
\centering
\includegraphics[height=2.6cm,width=8cm]{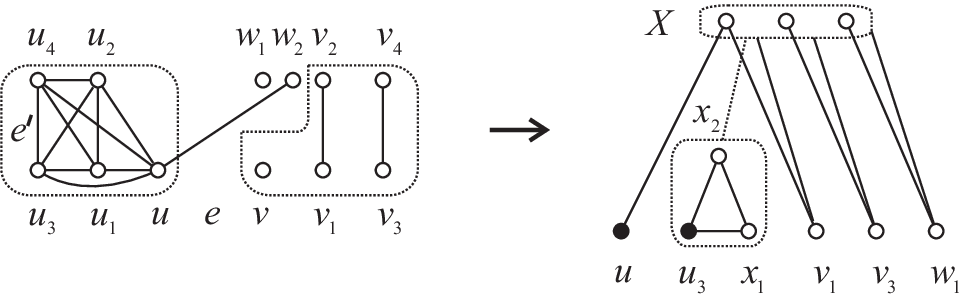}
\caption{\label{tu-7} $G_{uu_{3}}$ is $C13$ by applying Claim 1 to edge $uu_{3}$.}
\end{figure}

\smallskip
{\textbf{Case 6.} $G_{e}$ is $C7$.}

Let $M$ be a perfect matching of $G-S_{e}$ and $v_{1}v_{2}, v_{3}v_{4}\in M$.
We may assume that $ua_{1}\in E(G)$ and take an edge $e'=a_{1}v_{6}$.
Then $\overline{N[a_{1}]}\cap \overline{N[v_{6}]}\subseteq \{v, v_{1}, v_{2}, v_{3}, v_{4}\}$.
By Claim 2 and Cases 1 to 5, $G_{e'}$ is not $C1$, $C2$, $C3$, $C4$, $C5$, $C8$, $C9$,
$C10$, $C11$, $C12$, $C13$ or $C14$.

If $G_{e'}$ is $C6$, then $C6-X$ would contain an odd component $K_{3}$ induced
by three vertices in $\{v, v_{1}, v_{2}, v_{3}, v_{4}\}$, say $G[\{v_{1}, v_{2}, v_{3}\}]$.
Thus $C_{a_{1}}$ must be $G[\{a_{1}, u, v\}]$.
Since $G[\{v, v_{1}, v_{2}, v_{3}, v_{4}\}]$ is factor-critical,
$vv_{1}, vv_{2}$ or $vv_{3}\in E(G)$, a contradiction.

If $G_{e'}$ is $C7$, then we may assume $v_{1}, v_{3}$ (similarly, $v, v_{1}$)
as two trivial components of $C7-X$.
It follows that $a_{1}$ is the third trivial component of $C7-X$.
Otherwise, $v_{6}$ is the third trivial component and $v_{2}$ or $v_{4}\in$ $V(C_{a_{1}})$
but $v_{1}v_{2}, v_{3}v_{4}\in E(G)$, a contradiction.
Then $\{u, v_{5}\}\subseteq \overline{N[v_{1}]}= \{a_{1}, v_{3}\}\cup V(C_{v_{6}})$.
So $u, v_{5}\in V(C_{v_{6}})$ but $a_{1}u,$ $a_{1}v_{5}\in E(G)$, a contradiction.

\smallskip
{\textbf{Case 7.} $G_{e}$ is $C10$.}

Since $G-S_{e}$ has a perfect matching $M$,
assume that $v_{1}a_{1}, v_{6}a_{2}, v_{4}v_{5}\in M$. Consider edge $e'=v_{1}a_{2}$.
Clearly, $\overline{N[v_{1}]}\cap \overline{N[a_{2}]}\subseteq \{u, v, v_{2}, v_{3}, v_{4}, v_{5}\}$
and $uv, v_{2}v_{3}, v_{4}v_{5}$ are three independent edges. We apply Claim 1 to $e'$.
By Claim 2 and Cases 1 to 6, $G_{e'}$ is not $C1$, $C2$, $C3$, $C4$, $C5$, $C7$, $C8$,
$C9$, $C11$, $C12$ or $C14$.

If $G_{e'}$ is $C6$, then $G[\{v, v_{2}, v_{3}\}]$ must be an
odd component of $C6-X$. Assume that the others
are $G[\{a_{2}, v_{4}, v_{6}\}]$ and $G[\{v_{1}, x_{1}, x_{2}\}]$.
Hence $\overline{N[v_{2}]}\supseteq \{a_{2}, u, v_{1}, v_{4}, v_{5}, v_{6}, x_{1}, x_{2}\}$.
Since $v_{1}x_{1}$, $v_{1}x_{2}\in E(G)$ and $v_{1}u$, $v_{1}v_{5}\notin E(G)$,
$x_{1}, x_{2}\notin \{u, v_{5}\}$. So $d_{G}(v_{2})$$\leq n-9$, a contradiction.

If $G_{e'}$ is $C10$, similarly, $G[\{v, v_{2}, v_{3}\}]$ is an odd component of $C10-X$.
Assume that $v_{4}$ belongs to a trivial odd component.
It follows that $a_{2}$ (resp. $v_{1}$) belongs to the trivial (resp. nontrivial)
odd component of $C10-X$. Otherwise, $C_{a_{2}}$ would be $G[\{a_{2}, v_{5}, v_{6}\}]$
but $v_{4}v_{5}$, $v_{4}v_{6}\in E(G)$, a contradiction.
Let $C_{v_{1}}$ be $G[\{v_{1}, x_{1}, x_{2}\}]$.
Then $\overline{N[v_{4}]}\supseteq\{a_{2}, x_{1}, x_{2}, v, v_{1}, v_{2}, v_{3}, u\}$.
Since $v_{1}x_{1}$, $v_{1}x_{2}\in E(G)$ and $v_{1}u\notin E(G)$, $u\notin \{x_{1}, x_{2}\}$.
So $d_{G}(v_{4})$$\leq n-9$, a contradiction.

If $G_{e'}$ is $C13$, then we may assume that $\{u, v_{2}, v_{4}\}$
is an independent set of $G$, which induce three trivial odd components of $C13-X$.
It follows that $a_{2}$ (resp. $v_{1}$)
belongs to the trivial (resp. nontrivial) odd component of $C13-X$. Otherwise,
$C_{a_{2}}$ is either $G[\{a_{2}, v, v_{3}\}]$ or $G[\{a_{2}, v_{5}, v_{6}\}]$
but $uv$, $v_{4}v_{5}\in E(G)$, a contradiction.
Assume that $C_{v_{1}}$ is $G[\{v_{1}, x_{1}, x_{2}\}]$.
Then $\overline{N[v_{4}]}\supseteq\{u, v, v_{1}, v_{2}, v_{3}, a_{2}, x_{1}, x_{2}\}$.
Since $v_{1}x_{1}$, $v_{1}x_{2}\in E(G)$ and $v_{1}v$, $v_{1}v_{3}\notin E(G)$,
$x_{1}, x_{2}\notin \{v, v_{3}\}$. So $d_{G}(v_{4})$$\leq n-9$, a contradiction.

\smallskip
{\textbf{Case 8.} $G_{e}$ is $C6$.}

Let $M$ be a perfect matching of $G-S_{e}$. Assume that $v_{3}v_{4}, av_{5}\in M$.
We claim that $av_{3}, av_{4}\in E(G)$. Otherwise, say $av_{3}\notin E(G)$.
For edge $v_{3}v_{5}$, we have
$\overline{N[v_{3}]}\cap \overline{N[v_{5}]}=\{u, v, u_{1}, u_{2}, v_{1}, v_{2}\}$
and $uv, u_{1}u_{2}, v_{1}v_{2}$ are three independent edges.
We can give a similar discussion as Case 1 by using Claim 1 to $v_{3}v_{5}$.

If $au_{1}$ or $au_{2}\notin E(G)$, say $au_{1}\notin E(G)$, then
$\overline{N[u]}\cap \overline{N[u_{1}]}\subseteq \{a, v_{1}, v_{2}, v_{3}, v_{4}, v_{5}\}$.
Consider edge $uu_{1}$. By Claim 2 and Cases 1 to 7, $G-uu_{1}-S_{uu_{1}}$ may be $C6$, $C8$, $C11$, $C12$,
$C13$ or $C14$. Since $v_{1}v_{2}, v_{3}v_{4}, av_{5}$ are three independent edges,
$G-uu_{1}-S_{uu_{1}}$ is not $C8$, $C11$, $C12$ or $C14$.
Further, $G[\{a, v_{3}, v_{4}, v_{5}\}]$ is a $K_{4}$ and $v_{1}v_{2}\in E(G)$.
There is not an independent set with three vertices in $\{a, v_{1}, v_{2}, v_{3}, v_{4}, v_{5}\}$.
So $G-uu_{1}-S_{uu_{1}}$ is not $C13$. Thus $G-uu_{1}-S_{uu_{1}}$ would only be $C6$.
Then $C_{u}$ must be $G[\{u, v, a\}]$ and $G[\{v_{3}, v_{4}, v_{5}\}]$ would be an
odd component of $C6-X$. But $av_{3}, av_{4}, av_{5}\in E(G)$, a contradiction.
Thus $au_{1}$, $au_{2}\in E(G)$. Similarly, $av_{1}$, $av_{2}\in E(G)$.

Now consider edge $au_{1}$. Obviously, $\overline{N[a]}\cap \overline{N[u_{1}]}\subseteq \{v, w\}$,
where $w\in S_{e}$.  By Claim 2 and Cases 1, 3, 5 and 6,
$G-au_{1}-S_{au_{1}}$ is not any one of Configurations $C1$ to $C14$, which contradicts Claim 1.

\smallskip
{\textbf{Case 9.} $G_{e}$ is $C13$.}

For a perfect matching $M$ of $G-S_{e}$, we may assume that $b_{1}v_{3}, b_{2}v_{4}, b_{3}v_{5}\in M$.
Since $\delta(G-S_{e})\geq 2$, assume that $b_{1}u, b_{2}v_{3}\in E(G)$.
Take an edge $e'=b_{2}v_{3}$ and apply Claim 1 to $e'$.
We divide the proof into the two subcases as show in Fig. 12.

\begin{figure}[h]
\centering
\includegraphics[height=3.3cm,width=7cm]{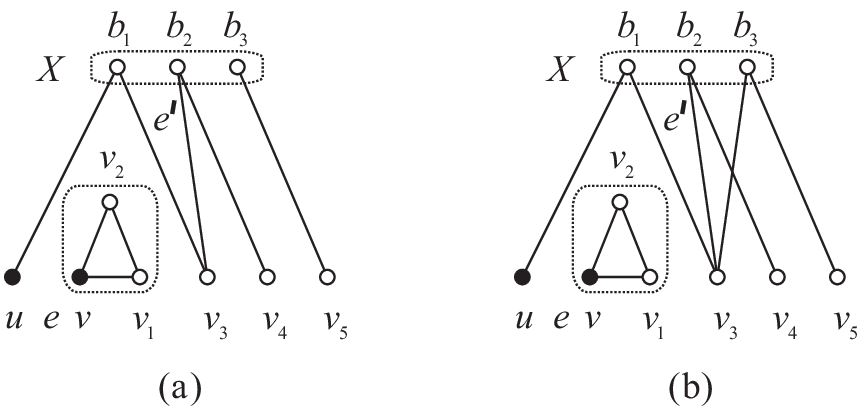}
\caption{\label{tu-05} The two configurations of $C13$.}
\end{figure}

{\textbf{Subcase 9.1.} $b_{3}v_{3}\notin E(G)$ (see Fig. 12 (a)).}

Then $\overline{N[b_{2}]}\cap \overline{N[v_{3}]}\subseteq \{u, v, v_{1}, v_{2}, v_{5}, b_{3}\}$
and $uv, v_{1}v_{2}, v_{5}b_{3}$ are three independent edges.
By Claim 2 and Cases 1 to 8, $G_{e'}$ would be $C13$. Assume that $\{u, v_{1}, v_{5}\}$ is
an independent set of $G$. It follows that $b_{2}$ is the fourth trivial odd component of $C13-X$.
Otherwise, $v_{3}$ is the fourth trivial odd component and
$C_{b_{2}}$ is either $G[\{b_{2}, v, v_{2}\}]$ or $G[\{b_{2}, b_{3}, v_{4}\}]$
but $uv, v_{5}b_{3}\in E(G)$, a contradiction.
Let $C_{v_{3}}$ be $G[\{v_{3}, x_{1}, x_{2}\}]$.
Then $\overline{N[u]}\supseteq \{v_{1}, v_{2}, v_{3}, v_{4}, v_{5}, x_{1}, x_{2}, b_{2}\}$.
Since $v_{3}x_{1}, v_{3}x_{2}\in E(G)$ and $v_{3}v_{2}, v_{3}v_{4}\notin E(G)$,
$x_{1}, x_{2}\notin \{v_{2}, v_{4}\}$. So $d_{G}(u)\leq n-9$, a contradiction.

{\textbf{Subcase 9.2.} $b_{3}v_{3}\in E(G)$ (see Fig. 12 (b)).}

Then $\overline{N[b_{2}]}\cap \overline{N[v_{3}]}\subseteq \{u, v, v_{1}, v_{2}, v_{5}, w\}$,
where $w\in S_{e}$.

If $b_{2}v_{5}\notin E(G)$, then $\overline{N[v_{5}]}=\{b_{2}, u, v, v_{1}, v_{2}, v_{3}, v_{4}\}$.
So $v_{5}w\in E(G)$ and $uv, v_{1}v_{2}, v_{5}w$ are three independent edges.
The proof is similar to Subcase 9.1.

If $b_{2}v_{5}\in E(G)$,
then $\overline{N[b_{2}]}\cap \overline{N[v_{3}]}\subseteq \{u, v, v_{1}, v_{2}, w\}$.
Since $uv, v_{1}v_{2}\in E(G)$, $G_{e'}$ would only be $C13$ by Claim 2 and Cases 1 to 8.
So we may assume that $\{u, v_{1}, w\}$ is an independent set of $G$.
Hence $b_{2}$ (resp. $v_{3}$) belongs to the trivial (resp. nontrivial) odd component of $C13-X$.
Otherwise, $C_{b_{2}}$ is $G[\{b_{2}, v, v_{2}\}]$ but $uv\in E(G)$, a contradiction.
Let $C_{v_{3}}$ be $G[\{v_{3}, x_{1}, x_{2}\}]$.
Then $\overline{N[u]}\supseteq \{v_{1}, v_{2}, v_{3}, v_{4}, v_{5}, x_{1}, x_{2}, b_{2}, w\}$.
Since $v_{3}x_{1}, v_{3}x_{2}\in E(G)$ and $v_{3}v_{2}, v_{3}v_{4}, v_{3}v_{5}\notin E(G)$,
$x_{1}, x_{2}\notin \{v_{2}, v_{4}, v_{5}\}$. So $d_{G}(u)\leq n-10$, a contradiction.

\smallskip
{\textbf{Case 10.} $G_{e}$ is $C12$.}

Let $M$ be a perfect matching of $G-S_{e}$.
Assume that $v_{1}v_{2}, b_{1}v_{3}, b_{2}v_{4}, b_{3}v_{5}\in M$.
Since $v_{4}$ has at least two neighbors in $X$, we discuss the three subcases as shown in Fig. 13.

\begin{figure}[h]
\centering
\includegraphics[height=3.2cm,width=9.6cm]{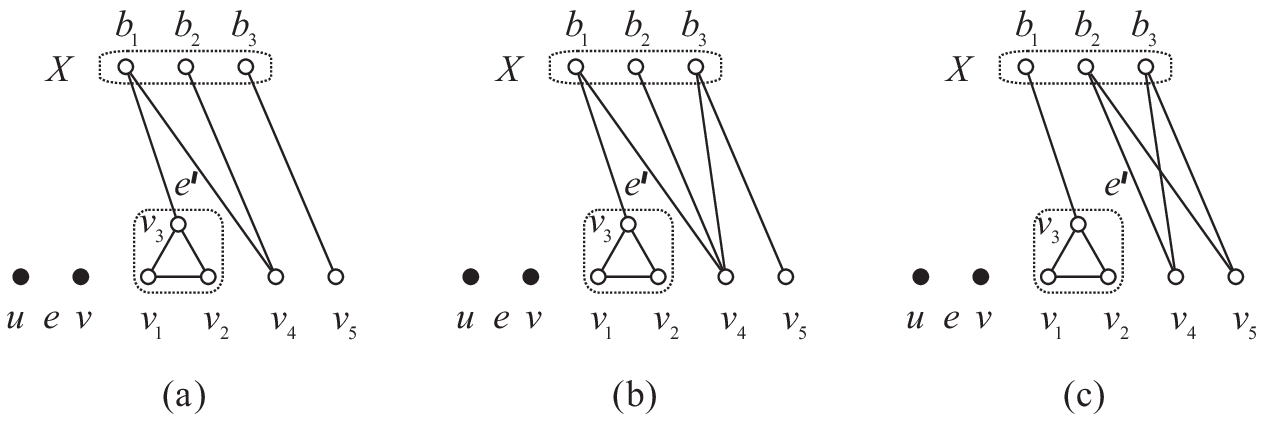}
\caption{\label{tu-06} The three configurations of $C12$.}
\end{figure}

{\textbf{Subcase 10.1.} $b_{1}v_{4}\in E(G)$, $b_{3}v_{4}\notin E(G)$ (see Fig. 13 (a)).}

Let $e'=b_{1}v_{4}$.
Then $\overline{N[b_{1}]}\cap \overline{N[v_{4}]}\subseteq$ $\{u, v, v_{1}, v_{2}, v_{5}, b_{3}\}$.
We apply Claim 1 to $e'$.
Since $uv, v_{1}v_{2}, v_{5}b_{3}$ are three independent edges, $G_{e'}$ is not $C8$, $C11$, $C12$ or $C14$.
By Claim 2 and Cases 1 to 9, we exclude the remaining configurations. This contradicts Claim 1.

{\textbf{Subcase 10.2.} $b_{1}v_{4}$, $b_{3}v_{4}\in E(G)$ (see Fig. 13 (b)).}

Consider edge $e'=b_{1}v_{4}$.
Then $\overline{N[b_{1}]}\cap \overline{N[v_{4}]}\subseteq \{u, v, v_{1}, v_{2}, v_{5}, w\}$,
where $w\in S_{e}$.

If $b_{1}v_{5}\notin E(G)$,
then $\overline{N[v_{5}]}=\{u, v, v_{1}, v_{2}, v_{3}, v_{4}, b_{1}\}$.
So $v_{5}w\in E(G)$ and $uv, v_{1}v_{2}, v_{5}w$ are three independent edges.
By similar discussions with Subcase 10.1, $G_{e'}$ is not any
one of Configurations $C1$ to $C14$, which contradicts Claim 1.

If $b_{1}v_{5}\in E(G)$,
then $\overline{N[b_{1}]}\cap \overline{N[v_{4}]}\subseteq \{u, v, v_{1}, v_{2}, w\}$.
Since $uv, v_{1}v_{2}$ are two independent edges, $G_{e'}$ is not $C12$ or $C14$.
By Claim 2 and Cases 1 to 9, $G_{e'}$ is not the remaining configurations,
which contradicts Claim 1.

{\textbf{Subcase 10.3.} $b_{1}v_{4}\notin E(G)$, $b_{3}v_{4}\in E(G)$ (see Fig. 13 (c)).}

Assume that $b_{1}v_{5}\notin E(G)$. Otherwise, we consider edge $b_{1}v_{5}$ the same as
Subcases 10.1 or 10.2. So $b_{2}v_{5}\in E(G)$.
Let $e'=b_{2}v_{4}$.
Then $\overline{N[b_{2}]}\cap \overline{N[v_{4}]}\subseteq \{b_{1}, u, v, v_{1}, v_{2}, v_{3}\}$
and $uv, v_{1}v_{2}, b_{1}v_{3}$ are three independent edges.
By Claim 2 and Cases 1 to 9, $G_{e'}$ is not any one of Configurations $C1$ to $C14$,
which is a contradiction to Claim 1.

\smallskip
{\textbf{Case 11.} $G_{e}$ is $C14$.}

For a perfect matching $M$ of $G-S_{e}$, we
assume that $c_{1}v_{1}, c_{2}v_{2}, c_{3}v_{3}, c_{4}v_{4}\in M$.
Since $\delta(G-S_{e})\geq 2$, assume that $c_{3}v_{4}\in E(G)$.
Let $e'=c_{3}v_{4}$. We apply Claim 1 to $e'$.
We divide the proof into the three subcases as shown in Fig. 14.

\begin{figure}[h]
\centering
\includegraphics[height=3.5cm,width=9.6cm]{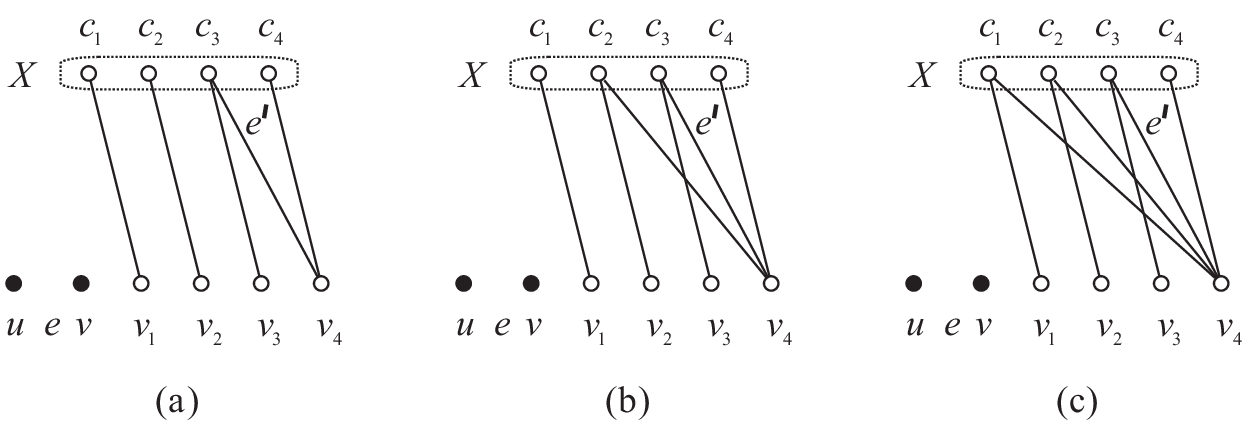}
\caption{\label{tu-07} The three configurations of $C14$.}
\end{figure}

{\textbf{Subcase 11.1.} $c_{1}v_{4}, c_{2}v_{4}\notin E(G)$ (see Fig. 14 (a)).}

It is easy to see that $\overline{N[c_{3}]}\cap \overline{N[v_{4}]}\subseteq \{c_{1}, c_{2}, u, v, v_{1}, v_{2}\}$
and $uv, v_{1}c_{1}, v_{2}c_{2}$ are three independent edges.
The proof is similar to Subcase 10.1.

{\textbf{Subcase 11.2.} $c_{1}v_{4}\notin E(G), c_{2}v_{4}\in E(G)$ (see Fig. 14 (b)).}

Then $\overline{N[c_{3}]}\cap \overline{N[v_{4}]}\subseteq \{c_{1}, u, v, v_{1}, v_{2}, w\}$,
where $w\in S_{e}$. By Claim 2 and Cases 1 to 10, $G_{e'}$ would be $C8$, $C11$ or $C14$.
Since $G[\{c_{1}, u, v, v_{1}, v_{2}, w\}]$ does not contain two disjoint triangles,
$G_{e'}$ is not $C11$.

If $G_{e'}$ is $C8$, then $v_{2}$ belongs to a trivial component of $C8-X$
and $G[\{c_{1}, u, v, v_{1}, w\}]$ is the nontrivial component
as $G[\{c_{1}, u, v, v_{1}, v_{2}\}]$ is not factor-critical.
Thus $\overline{N[v_{2}]}\supseteq\{c_{1},$ $c_{3}, u, v, v_{1}, v_{3}, v_{4}, w\}$.
So $d_{G}(v_{2})\leq n-9$, a contradiction.

If $G_{e'}$ is $C14$, we assume that $\{u, v_{1}, v_{2}, w\}$ is an independent set of $G$.
Then $\overline{N[v_{1}]}=\{c_{3}, u, v, v_{2}, v_{3}, v_{4}, w\}$.
So $v_{1}c_{2}$, $v_{1}c_{4}\in E(G)$.
Consider edge $c_{2}v_{4}$ and apply Claim 1 to edge $c_{2}v_{4}$.
Then $\overline{N[c_{2}]}\cap \overline{N[v_{4}]}\subseteq \{c_{1}, u, v, v_{3}, w\}$.
By Claim 2 and Cases 1 to 10, $G-c_{2}v_{4}-S_{c_{2}v_{4}}$ is also $C14$.
Let $\{u, c_{1}, v_{3}, w\}$ be an independent set of $G$.
Thus $\overline{N[w]}\supseteq \{c_{1}, c_{2}, c_{3}, u, v_{1}, v_{2}, v_{3}, v_{4}\}$.
So $d_{G}(w)$$\leq n-9$, a contradiction.

{\textbf{Subcase 11.3.} $c_{1}v_{4}, c_{2}v_{4}\in E(G)$ (see Fig. 14 (c)).}

Then $\overline{N[c_{3}]}\cap \overline{N[v_{4}]}\subseteq \{u, v, v_{1}, v_{2}, w_{1}, w_{2}\}$,
where $w_{1}, w_{2}\in S_{e}$.
By Claim 2 and Cases 1 to 10, $G_{e'}$ would only be $C14$.
So we need to find an independent set $T$ with size four in $\{u, v, v_{1}, v_{2}, w_{1}, w_{2}\}$.
Since $uv\in E(G)$, it is impossible that $w_{1}, w_{2} \notin T$.
We claim that one of $\{w_{1}, w_{2}\}$ belongs to $T$. Otherwise,
if $w_{1}, w_{2} \in T$, then $v_{1}$ or $v_{2}\in T$, say $v_{1}\in T$.
Thus $\overline{N[v_{1}]}\supseteq \{u, v, v_{2}, v_{3}, v_{4}, c_{3}, w_{1}, w_{2}\}$.
So $d_{G}(v_{1})\leq n-9$, a contradiction.
Assume that $w_{1}\in T$ and $w_{2}\notin T$. So $v_{1}$, $v_{2}\in T$.
Then we may assume that $T=\{u, v_{1}, v_{2}, w_{1}\}$.

Since $\overline{N[v_{1}]}=\{u, v, v_{2}, v_{3}, v_{4}, c_{3}, w_{1}\}$,
$v_{1}c_{2}$, $v_{1}c_{4}\in E(G)$. For edge $c_{2}v_{4}$, we have
$\overline{N[c_{2}]}\cap \overline{N[v_{4}]}$ $\subseteq \{u, v, v_{3}, w_{1}, w_{2}\}$.
Then $G-c_{2}v_{4}-S_{c_{2}v_{4}}$ is still $C14$ by using Claim 1 to edge $c_{2}v_{4}$.
Let $\{u, v_{3}, w_{1}, w_{2}\}$ be an independent set of $G$.
Thus $\overline{N[v_{3}]}\supseteq \{u, v, v_{1}, v_{2}$, $v_{4}, w_{1}, w_{2}, c_{2}\}$.
So $d_{G}(v_{3})\leq n-9$, a contradiction.

\smallskip
{\textbf{Case 12.} $G_{e}$ is $C8$.}

Assume that $a_{1}v_{1}, a_{2}v_{2}, v_{3}v_{4}, v_{5}v_{6}$
belong to a perfect matching of $G-S_{e}$. Let $e'=v_{1}a_{2}$.
Then $\overline{N[v_{1}]}\cap$ $\overline{N[a_{2}]}\subseteq \{u, v, v_{3}, v_{4}, v_{5}, v_{6}\}$.
We apply Claim 1 to $e'$.
Since $uv, v_{3}v_{4}, v_{5}v_{6}$ are three independent edges, $G_{e'}$ is not $C8$ or $C11$.
By Claim 2 and Cases 1 to 11, $G_{e'}$ is not the remaining configurations.
This is a contradiction to Claim 1.

\smallskip
{\textbf{Case 13.} $G_{e}$ is $C11$.}

Let $a_{1}v_{3}, a_{2}v_{6}, v_{1}v_{2}, v_{4}v_{5}$ belong to a perfect matching
of $G-S_{e}$. We may assume that $ua_{1}\in E(G)$. Let $e'=ua_{1}$. We apply Claim 1 to $e'$.

If $ua_{2}\notin E(G)$,
then $\overline{N[u]}\cap \overline{N[a_{1}]}\subseteq \{v_{1}, v_{2}, v_{4}, v_{5}, v_{6}, a_{2}\}$
and $v_{1}v_{2}, v_{4}v_{5}, a_{2}v_{6}$ are three independent edges.
It is obvious that $G_{e'}$ is not any one of Configurations $C1$ to $C14$,
which contradicts Claim 1.

If $ua_{2}\in E(G)$,
then $\overline{N[u]}\cap \overline{N[a_{1}]}\subseteq \{v_{1}, v_{2}, v_{4}, v_{5}, v_{6}, w\}$,
where $w\in S_{e}$. By Claim 2 and Cases 1 to 12, $G_{e'}$ would only be $C11$.
Then the two nontrivial odd components of $C11-X$ must be
$G[\{v_{1}, v_{2}, w\}]$ and $G[\{v_{4}, v_{5}, v_{6}\}]$.
Thus $\overline{N[v_{4}]}=\{a_{1}, u, v, v_{1}, v_{2}, v_{3}, w\}$. So $a_{2}v_{4}\in E(G)$.
Similarly, $a_{2}v_{5}\in E(G)$. Now take another edge $ua_{2}$ and apply Claim 1 to $ua_{2}$.
Then $\overline{N[u]}\cap \overline{N[a_{2}]}\subseteq \{v_{1}, v_{2}, v_{3}, w\}$.
It is easy to see that $G-ua_{2}-S_{ua_{2}}$ is not any one of Configurations $C1$ to $C14$,
which contradicts Claim 1.

\smallskip
{\textbf{Case 14.} $G_{e}$ is $C3$.}

Assume that $av_{1}, v_{2}v_{3}, v_{4}v_{5}, v_{6}v_{7}$ belong to a perfect matching
of $G-S_{e}$. Let $e'=ua$.
Then $\overline{N[u]}\cap \overline{N[a]}\subseteq \{v_{2}, v_{3}, v_{4}, v_{5}, v_{6}, v_{7}\}$.
We apply Claim 1 to $e'$. It is easy to see that $G_{e'}$ is not $C3$.
By Cases 1 to 13, $G_{e'}$ is not the other configurations.
This is a contradiction to Claim 1.

\smallskip
Combining Cases 1 to 14, we complete the proof.
~~~~~~~~~~~~~~~~~~~~~~~~~~~~~~~~~~~~~~~~~~~~~~~~~~~~~~~~~~~~$\square$

\end{document}